\documentclass{article}
\usepackage{graphicx}
\usepackage{amsmath}
\usepackage{amssymb}
\usepackage{amsthm}
\usepackage{mathrsfs}
\newtheorem{dfn}{Definition}[section]
\newtheorem{thm}[dfn]{Theorem}
\newtheorem{prop}[dfn]{Proposition}

\newtheorem{rem}[dfn]{Remark}

\numberwithin{equation}{section}

\begin{document}

\title{Numerical validation of blow-up solutions of ordinary differential equations}

\author{Akitoshi Takayasu 
\thanks{Faculty of Engineering, Information and Systems, University of Tsukuba, 1-1-1 Tennodai, Tsukuba, Ibaraki 305-8573, Japan ({\tt takitoshi@aoni.waseda.jp})},
Kaname Matsue
 \thanks{The Institute of Statistical Mathematics, 10-3 Midori-cho, Tachikawa, Tokyo 190-8562, Japan},
Takiko Sasaki
 \thanks{Global Education Center, Waseda University, 1-104 Totsuka-machi, Shinjuku, Tokyo 169-8050, Japan},
Kazuaki Tanaka
\thanks{Graduate School of Fundamental Science and Engineering, Waseda University, 3-4-1 Okubo, Shinjuku, Tokyo 169-8555, Japan},\\
Makoto Mizuguchi
\thanks{Graduate School of Fundamental Science and Engineering, Waseda University, 3-4-1 Okubo, Shinjuku, Tokyo 169-8555, Japan},\ 
and
Shin'ichi Oishi
\thanks{Department of Applied Mathematics, Faculty of Science and Engineering, Waseda University, 3-4-1 Okubo, Shinjuku, Tokyo 169-8555, Japan}
}
\maketitle

\begin{abstract}
This paper focuses on blow-up solutions of ordinary differential equations (ODEs).
We present a method for validating blow-up solutions and their blow-up times, which is based on compactifications and the Lyapunov function validation method.
The necessary criteria for this construction can be verified using 
interval arithmetic techniques.
Some numerical examples are presented to demonstrate the applicability of our method.
\end{abstract}


\section{Introduction}
In this paper, we consider the initial value problem defined by the following ordinary differential equations in $\mathbb{R}^m$ ($m\in\mathbb{N}$):
\begin{equation}\label{eqn:ODE}
	\frac{dy(t)}{dt}=f\left(y(t)\right),~y(0)=y_0,
\end{equation}
where $t\in[0,T)$ with $0<T\le\infty$, $f:\mathbb{R}^m\to\mathbb{R}^m$ is a $C^1$ function, and $y_0\in\mathbb{R}^m$.
Unless otherwise noted, $f$ is assumed to be a polynomial, whose coefficients are real numbers.
Our focus in this paper is a class of solutions of (\ref{eqn:ODE}) called {\em blow-up solutions}.
\begin{dfn}\label{def:blowup}\rm
Define $t_{\max}>0$ as
\[
	t_{\max}:=\sup\left\{\bar t:\mbox{a solution $y\in C^1([0,\bar t))$ of \eqref{eqn:ODE} exists}\right\}.
\]
We say that the solution $y$ of \eqref{eqn:ODE} {\em blows up} if $t_{\max}<\infty$.
In such a case, $t_{\max}$ is called the {\em blow-up time} of \eqref{eqn:ODE}.
\end{dfn}
The simplest example of blow-up phenomena can be seen for the following ordinary differential equation (ODE) in $\mathbb{R}^1$:
\begin{equation}
\label{simplest}
\frac{dy}{dt} = y^2,\quad y(0) = y_0.
\end{equation}
When $y_0 > 0$, the exact solution of (\ref{simplest}) is $y(t) = (y_0^{-1}-t)^{-1}$. 
The value of $y(t)$ becomes infinite as $t\to y_0^{-1}-0$.
That is, $y(t)$ blows up at $t = y_0^{-1}$.
Blow-up solutions can also be observed for partial differential equations (PDEs), such as the nonlinear heat equations (e.g., \cite{bib:F1966, bib:K1996}) given by
\begin{equation}
\label{eqn:heat_intro}
u_t = \Delta u + |u|^{p-1}u,\quad p > 1,
\end{equation}
the nonlinear wave equations (e.g., \cite{bib:J1981, bib:L1974}), and the nonlinear Schr\"{o}dinger equations (e.g., \cite{bib:G1977}). 
In the case of PDEs, many researchers have studied blow-up phenomena such as blow-up times, blow-up criteria, the behavior of solutions near blow-up times (e.g., blow-up rate), and the topology or geometry of blow-up sets.
Studies of blow-up phenomena can be of importance both mathematically and physically.
For example, in the case of the nonlinear heat equation (\ref{eqn:heat_intro}), blow-up solutions describe the combustion of solid fuels \cite{bib:BE1989}. 
Similarly, the blow-up time corresponds to the time when the fuel ignites. 
Blow-up phenomena associated with (\ref{eqn:heat_intro}) thus describe the process of combustion.

\begin{figure}[htbp]
\begin{minipage}{1.0\hsize}
\centering
\includegraphics[width=6.0cm]{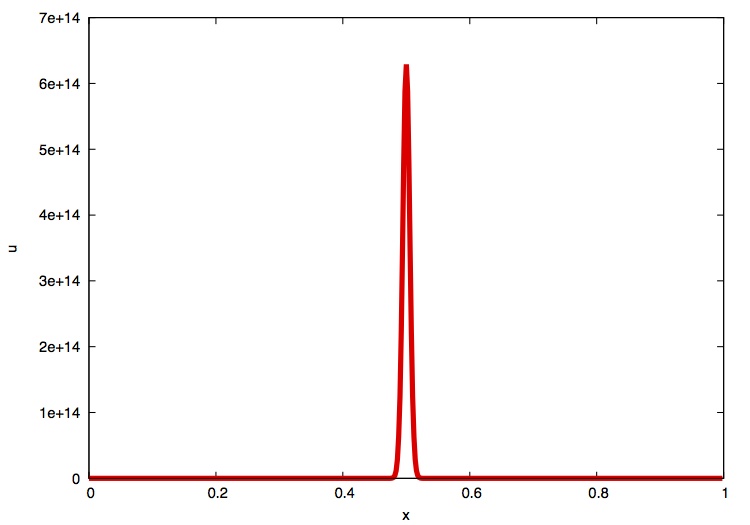}
\end{minipage}
\caption{Numerical blow-up solution of \eqref{eqn:heat_intro} with $p=3$.}
\label{fig:blow-up}
The $L^\infty$-norm of solutions becomes larger and larger, and may become infinite within a finite time.
We typically regard solutions whose $L^\infty$-norms become sufficiently large within finite times as ``blow-up solutions'' in a numerical sense.
\end{figure}

The numerical analysis of blow-up solutions, such as of nonlinear heat and reaction-diffusion equations \cite{bib:ALM998, bib:C1986, bib:CHO2007, bib:N1975, bib:NU1977}, of nonlinear wave equations \cite{bib:C2010, bib:SS2015},
and of nonlinear Schr\"{o}dinger equations \cite{bib:ADKM2003, bib:BCMS2008}, has also been studied.
However, in almost all numerical studies concerning blow-up solutions, ``blow-up solutions'' have been only computed approximately. 
For example, typical numerical computations of blow-up solutions begin by setting an appropriately large number $M$, say $10^6$. 
Then, one numerically solves the differential equations, and regards computed solutions whose supremum norms are larger than $M$ as blow-up solutions (Fig. \ref{fig:blow-up}). 
However, this criterion provides us with no proof that these computed blow-up solutions are {\em rigorous} blow-up solutions.
In other words, it is possible that ``numerical blow-up solutions'' just describe extremely large but {\em bounded} solutions. For example, consider (\ref{simplest}) again, and the perturbed equation
\begin{equation*}
\frac{dy}{dt} = y^2 - \epsilon y^3,\quad y(0) = y_0 > 0,
\end{equation*}
where $\epsilon > 0$ is a sufficiently small parameter.
One can easily see that the solution tends to $y = 1/\epsilon$ as $t\to \infty$.
Obviously, this solution is not a blow-up solution, while the dominant behavior of this solution resembles that of $dy/dt = y^2$.
In such a case for a general system, it is not easy to judge whether a computed solution is truly a blow-up solution.
Therefore, an exact criterion for blow-up solutions is necessary to concretely obtain rigorous blow-up solutions.

The blow-up time $t_{\max}$ is one of the key considerations for blow-up solutions.
Some specific solutions, such as self-similar solutions, can be described via transformations involving $t_{\max}$ (see, e.g., \cite{bib:FM2002}), in which case we assume that $t_{\max}$ is known in advance.
However, the detection of $t_{\max}$ itself is not easy, because $t_{\max}$, in general, depends on an initial condition (as can be observed in (\ref{simplest})), and the behavior of solutions near to blow-up times changes rapidly.
Although there are many numerical studies regarding the detection of blow-up times with sufficient accuracy (see, e.g., \cite{bib:CHO2007}), the problem of whether the computed solution really blows up still remains.
Thus,  for further analysis of blow-up phenomena from the numerical viewpoint, it would be helpful to have a procedure for detecting {\em rigorous} blow-up solutions and their blow-up times, or their explicit enclosures. 

\bigskip
Considering this background regarding blow-up solutions, 
we provide a method for numerically validating blow-up solutions of \eqref{eqn:ODE}, as well as their blow-up times.
The key point in our method is the way of interpreting blow-up solutions.
Blow-up solutions are not compatible with numerical computations without additional considerations, because their norm goes to infinity within a finite time.
To overcome this difficulty, we analyze blow-up solutions using two instruments. 
The first is the {\em compactification} of base spaces, following Poincar\'e (\cite{bib:P}. See also \cite{bib:Elias-Gingold, bib:H}).
The compactification embeds the original space into compact manifolds with boundaries, which enables us to describe ``infinity'' as points on the boundaries.
Elias and Gingold \cite{bib:Elias-Gingold} discuss an appropriate class of such a compactification so that dynamical systems on the compactified space including the boundary make sense.
Blow-up solutions are then regarded as the trajectories that tend to the boundary of the compactified space.
The ``infinity'' where a solution diverges in the original system is translated into a ``critical point at infinity'' in the transformed system.
Thus, we can consider a divergent solution as a global solution that is asymptotic to a critical point at infinity.
This global solution can be validated using standard approaches in the theory of dynamical systems (see, e.g., \cite{bib:Mat2, bib:MTKO2015, bib:W}).
The second instrument we employ is the validation of {\em Lyapunov functions}.
A Lyapunov function, say $L$, guarantees the monotonous behavior of dynamical systems in the domain of $L$.
The monotonicity of $L$ enables us to derive a re-parameterization of trajectories in terms of $L$ along solution curves, which is called {\em Lyapunov tracing} in \cite{bib:MHY}.
This re-parameterization around critical points at infinity yields explicit estimates of blow-up times.
Preceding works \cite{bib:Mat, bib:MHY} indicate that such explicit and rigorous estimates of blow-up times can be derived with computer assistance.
The above two instruments, together with standard methodologies of numerical validation and dynamical systems, lead to the validation of blow-up solutions along with their blow-up times.
\par
The applicability of our method is exhibited with computer assistance based on interval arithmetic techniques.
Interval arithmetic enables us to compute enclosures in which mathematically correct objects are contained.
In dynamical systems theory, there are many applications for interval arithmetic, including validations of global trajectories, determining stability of invariant sets, and determining parameter ranges where dynamical bifurcations occur (see, e.g., \cite{bib:AK2010, bib:CGL2015, bib:CZ2015, bib:Mat2, bib:Mat, bib:MHY, bib:MTKO2015, bib:W, bib:Zgl2010}).
In our case, {\em affine arithmetic} \cite{bib:kashikashi, bib:kv}, an enhanced version of interval arithmetic, is applied to validate explicit enclosures of blow-up solutions and their blow-up times.

\bigskip
This paper is organized as follows.
In Section \ref{section-compact}, we introduce a notion of addmissible compactifications.
In addition, we introduce the normalization of the time scale $t$ in \eqref{eqn:ODE}, so that dynamical systems including points at infinity are valid.
The vector field with normalized time scale is central to our considerations.
In Section \ref{section-validation}, we present a criterion for validating blow-up solutions. 
Here, the (locally defined) Lyapunov function is also introduced.
By verifying this criterion numerically,
we can validate blow-up solutions and explicit bounds of their blow-up times.
Some sample numerical validation results are presented in Section \ref{section-numerical} to demonstrate the applicability of our method. 
\section{Compactification of dynamical systems}
\label{section-compact}
The compactification of spaces, which we apply throughout this paper, is a technique that enables us to describe points ``at infinity'' explicitly.
Roughly speaking, compactification is an operation that embeds the original vector space homeomorphically into a compact manifold (possibly with the boundary).
The simplest example is the one-point compactification of $\mathbb{R}^m$ into the unit $m$-sphere $S^m$ (the Bendixson compactification \cite{bib:H}), in which case the north pole of $S^m$ corresponds to infinity for $\mathbb{R}^m$.
Our basis for validating blow-up solutions is the application of compactifications to dynamical systems.
Although there are several possible compactifications, we want to choose a ``good'' compactification.
In other words, we require a compactification in which a ``dynamical system at infinity'' is nondegenerate in an appropriate sense.
Elias and Gingold \cite{bib:Elias-Gingold} discuss a general class of compactifications that satisfy our requirement.

Throughout the remainder of this section, we will review this class of compactifications, as well as several examples.
Throughout the remainder of this paper, let $\langle \cdot, \cdot \rangle$ denote the standard inner product on $\mathbb{R}^m$, and let $\|\cdot \|$ be its associated norm; that is, the Euclidean norm.
Moreover, the nonlinearity $f$ is assumed to be a polynomial.
\subsection{General admissible compactifications}
There are many approaches to compactifying the total space $\mathbb{R}^m$. 
If we want to apply this idea to dynamical systems at infinity, then a natural aim is to reduce degeneracy as far as possible. 
In particular, we want to explicitly distinguish directions towards infinity.
This is achieved by using a mapping between $\mathbb{R}^m$ and the open unit ball $\mathcal{D}\subset \mathbb{R}^m$, and an identification of unit vectors on $\overline{\mathcal{D}}$ with directions at infinity. 
This mapping allows the attachment of a set of ``points at infinity" in any direction $p\in \partial \mathcal{D}$ to the whole space $\mathbb{R}^m$.
\begin{dfn}[\cite{bib:Elias-Gingold}]\rm
We say that a sequence of points $\{y_k\}_{k\geq 1}\subset \mathbb{R}^m$ {\em tends to infinity in the direction $y_\ast$, where $\|y_\ast\|=1$}, if it holds that $\|y_k\| \to \infty$ and $y_k / \|y_k\|\to y_\ast$ as $k\to \infty$.
\end{dfn}
 
Let $T : \mathbb{R}^m \to \mathcal{D}$ be a map given by
\begin{equation}
\label{compact-abstract}
T(y) = x := \frac{y}{\kappa(y)},\quad \kappa (y) = \kappa(y_1,\cdots, y_m),
\end{equation}
where $\kappa :\mathbb{R}^m\to \mathbb{R}_{>0}:=\{a\in\mathbb{R}:a>0\}$ is continuous. Equivalently,
\begin{equation}
\label{compact-component}
x_i := \frac{y_i}{\kappa(y_1,\cdots, y_m)}.
\end{equation}
We review a special class of such maps that preserve information regarding the original dynamical system on $\mathbb{R}^m$. The following admissibility conditions ensure a correspondence between the original dynamical system and the transformed one.
\begin{dfn}[\cite{bib:Elias-Gingold}]\rm
\label{dfn-adm}
The map $T$ given by (\ref{compact-abstract}) is called an {\em admissible compactification} if the following statements hold:
\begin{description}
\item[(A0) ] $\kappa(y) > \|y\|$.
\item[(A1) ] $\kappa(y) = O(\|y\|)$ as $\|y\|\to \infty$.
\item[(A2) ] $\nabla \kappa(y)\sim y/\|y\|$ as $\|y\|\to \infty$.
\item[(A3) ] $\langle y,\nabla \kappa(y) \rangle < \kappa(y)$.
\end{description}
\end{dfn}
The validity of Definition \ref{dfn-adm} is stated in \cite{bib:Elias-Gingold}.
In particular, $T$ is bijective.
Note that for explicit expressions of $\kappa$ it is often the case that $T$ becomes a homeomorphism.

\bigskip
The next step is to transform the differential equation (\ref{eqn:ODE}) via (\ref{compact-abstract}).
A direct application of (\ref{compact-abstract}) to (\ref{eqn:ODE}) yields
\begin{align*}
\frac{dx}{dt} &= \frac{d}{dt}(y/\kappa(y)) = \kappa^{-1}\frac{dy}{dt} - \kappa^{-2} \left\langle \nabla \kappa(y), \frac{dy}{dt}\right\rangle y\\
	&= \kappa^{-1}[f(y) - \kappa^{-1}\langle \nabla \kappa, f(y)\rangle y]\\
	&= \kappa^{-1}[f(\kappa x) - \langle \nabla \kappa , f(\kappa x)\rangle x],
\end{align*}
i.e.,
\begin{equation}
\label{transformed}
\frac{dx}{dt} = {[\kappa(T^{-1}(x))]^{-1}}[f(\kappa x) - \langle \nabla \kappa , f(\kappa x)\rangle x],~
x(0)=\frac{y_0}{\kappa(y_0)}.
\end{equation}

A correspondence of invariant sets between (\ref{eqn:ODE}) and (\ref{transformed}) is given as follows.
\begin{prop}[\cite{bib:Elias-Gingold}]
\label{prop-fixpt}
The transformation $T$ in \eqref{compact-abstract} maps bounded equilibria of \eqref{eqn:ODE} in $\mathbb{R}^m$ into equilibria of \eqref{transformed} in $\mathcal{D}$, and vice versa.
\end{prop}
Next, we focus on points at infinity. These are realized via $T$ by points on $\partial \mathcal{D}$. 
On $\partial \mathcal{D}$, $\kappa(T^{-1}(x))\to \infty$, and hence (\ref{transformed}) is obviously singular on $\partial \mathcal{D}$.
Now, from (\ref{transformed}) we extract its singular part and its continuous part.
Let $d$ be the degree of the polynomial $f(y)$, and write $f(y)$ as follows:
\begin{equation*}
f(y) = p_0(y) + p_1(y) + \cdots + p_d(y),
\end{equation*}
where $p_j(y), j=0,\cdots, d$ are homogeneous polynomials of degree $j$, respectively.
Regarding the behavior of each $p_i(y)$ for large $\|y\|$, we expect that the highest order term $p_d(y)$ governs the behavior of $f(y)$. 
With this observation in mind, we set $y=\kappa x$, and define
\begin{equation}
\label{tilde-f}
\tilde f(x,\kappa) := \kappa^{-d}f(\kappa x) = \kappa^{-d}p_0(x) + \kappa^{-d+1}p_1(x) + \cdots + p_d(x),
\end{equation}
where $\kappa = \kappa(T^{-1}(x))$. According to (\ref{tilde-f}), (\ref{transformed}) can be rewritten as
\begin{equation}
\label{normalized-aux}
\frac{dx}{dt} = {[\kappa(T^{-1}(x))]^{d-1}}[\tilde f(x, \kappa) - \langle \nabla \kappa , \tilde f(x, \kappa)\rangle x],~
x(0)=\frac{y_0}{\kappa(y_0)}.
\end{equation}
$\kappa^{d-1}$ is positive and unbounded on $\partial \mathcal{D}$, while the rest of (\ref{normalized-aux}) is continuous on the {\em closed} ball $\overline{\mathcal{D}}$. If we define a new independent time variable $\tau$ along a trajectory $y(t)$, or equivalently $x(t)$, by
\begin{equation}
\label{time-transform}
\frac{d\tau}{dt} = \kappa(y(t))^{d-1},
\end{equation}
i.e.,
\begin{equation*}
\tau = \int_0^t (\kappa(y(s)))^{d-1}ds,
\end{equation*}
then we can write (\ref{normalized-aux}) uisng the $\tau$-time variable as
\begin{equation}
\label{normalized}
\frac{dx}{d\tau} = \tilde f(x, \kappa) - \langle \nabla \kappa , \tilde f(x, \kappa)\rangle x \equiv g(x),~
x(0)=\frac{y_0}{\kappa(y_0)}.
\end{equation}
Obviously, the vector field $g$ is defined on $\overline{\mathcal{D}}$, not only on $\mathcal{D}$.
Trajectories of (\ref{eqn:ODE}) in $\mathbb{R}^m$ and those of (\ref{normalized}) in $\mathcal{D}$ have the same topology and, according to Proposition \ref{prop-fixpt}, the same (bounded) equilibria. 
In order to unite bounded equilibria of (\ref{eqn:ODE}) and ``equilibria'' at infinity in the same framework, the following notion is proposed.
\begin{dfn}[Critical point at infinity. See, e.g., \cite{bib:Elias-Gingold, bib:H}]\rm
We say that (\ref{eqn:ODE}) has {\em a critical point at infinity in the direction $x_\ast$, $\|x_\ast\|=1$}, if $x_\ast$ is an equilibrium of (\ref{normalized}) in $\partial \mathcal{D}$; namely, the right-hand side of (\ref{normalized}) vanishes there.
\end{dfn}
\begin{rem}
The terminology ``critical point'' refers to an equilibrium point in a dynamical system, which may cause a little confusion for researchers of dynamical systems. 
Following the convention used in preceding works, we shall only use the terminology ``critical point'' to refer to critical points at infinity.
\end{rem}
The following proposition indicates that solutions of (\ref{eqn:ODE}) that diverge can be described in terms of critical points at infinity.
\begin{prop}[\cite{bib:Elias-Gingold}]
\label{prop-asymptotic}
Assume that a solution $y(t)$ of \eqref{eqn:ODE} has a maximal interval of existence $(a,b)\subset \mathbb{R}$, and that $y(t)$ tends to infinity in the direction $x_\ast$ as $t\to b-0$ (or as $t\to a+0$).
Then, $x_\ast$ is an equilibrium of \eqref{normalized} on $\partial \mathcal{D}$. 
\end{prop}
Notice that $a$ and $b$ admit $-\infty$ and $+\infty$, respectively, which implies that this property itself is not sufficient to describe blow-up solutions. 
For a precise distinction, we shall call a solution $y(t)$ that diverges as $t\to \infty$, i.e., $b = +\infty$ in Proposition \ref{prop-asymptotic}, a {\em grow-up solution}.

\bigskip
Observe that $\|x\|\to 1$ corresponds to the situation that $\|y\| \to \infty$. 
By assumptions (A1) and (A2), we know that $\nabla \kappa(y) \sim y/\|y\| \sim y/\kappa(y) = x$.
Thus, the equation $g(x) = 0$ (in (\ref{normalized})) with $\|x\|=1$ is equivalent to
\begin{equation}
\label{eq-critical}
\tilde f(x, \infty) - \langle x , \tilde f(x, \infty)\rangle x = 0.
\end{equation}
We shall refer to (\ref{eq-critical}) as the {\em (nonlinear) characteristic equation of (\ref{eqn:ODE}) at infinity}. 
Note that it holds that $\kappa\to \infty$ as $x$ goes to $\partial \mathcal{D}$, hence,
\begin{equation}
\tilde f(x,\kappa)\to p_d(x_\ast)
\quad
\text{ as }
\quad
x\to x_\ast\in \partial \mathcal{D}.
\end{equation}
As a consequence, critical points at infinity can be determined by just the highest order term $p_d(x)$.
This fact is summarized with additional observations by the following proposition.
\begin{prop}[\cite{bib:Elias-Gingold}]
\begin{enumerate}
\item The definition of critical points at infinity of \eqref{eqn:ODE} is independent of the choice of admissible compactification \eqref{compact-abstract}.
\item A critical point at infinity in the $x$ direction depends only on the highest order polynomial $p_d$ of the original vector field $f$. The characteristic equation at infinity \eqref{eq-critical} is equivalent to 
\begin{equation*}
p_d(x) - \langle x,p_d(x)\rangle x = 0,\quad \|x\|=1.
\end{equation*}
\item If \eqref{eqn:ODE} has a critical point at infinity in the direction $x_\ast$, then $-x_\ast$ is also a critical point at infinity of \eqref{eqn:ODE}.
\end{enumerate}
\end{prop}
A series of facts concerning compactifications leads to the translation of blow-up solutions (and grow-up solutions) in terms of dynamical systems as global solutions that are asymptotic to 
equilibria of (\ref{normalized}) on $\partial \mathcal{D}$.
This translation presents the possibility of validating these solutions with computer assistance, as in previous studies such as \cite{bib:Mat2, bib:MTKO2015, bib:W}. 
However, note that the above results do not yet distinguish blow-up solutions from grow-up solutions.
If we want to detect blow-up solutions, we need an additional procedure for computing blow-up times, which is discussed in Section \ref{section-validation}.
 
\subsection{Poincar\'{e} compactification}
\label{section-poincare}
In practical computations, we choose an admissible compactification explicitly.
The Poincar\'{e} compactification is one of the simplest admissible compactifications.
\begin{dfn}\rm
For the unit open ball $\mathcal{D}\subset \mathbb{R}^m$, the {\em Poincar\'{e} compactification of $\mathbb{R}^m$} is given by the map
\begin{equation}
\label{poincare}
T_{Poin}(y) = x := \frac{y}{\sqrt{1+\|y\|^2}},\quad \text{i.e. }\kappa(y_1,\cdots, y_m) = \sqrt{1+\|y\|^2} = \sqrt{1+\sum_{i=1}^m y_i^2}.
\end{equation}
\end{dfn}
\begin{figure}[h]
\centering
\includegraphics[width=6.0cm]{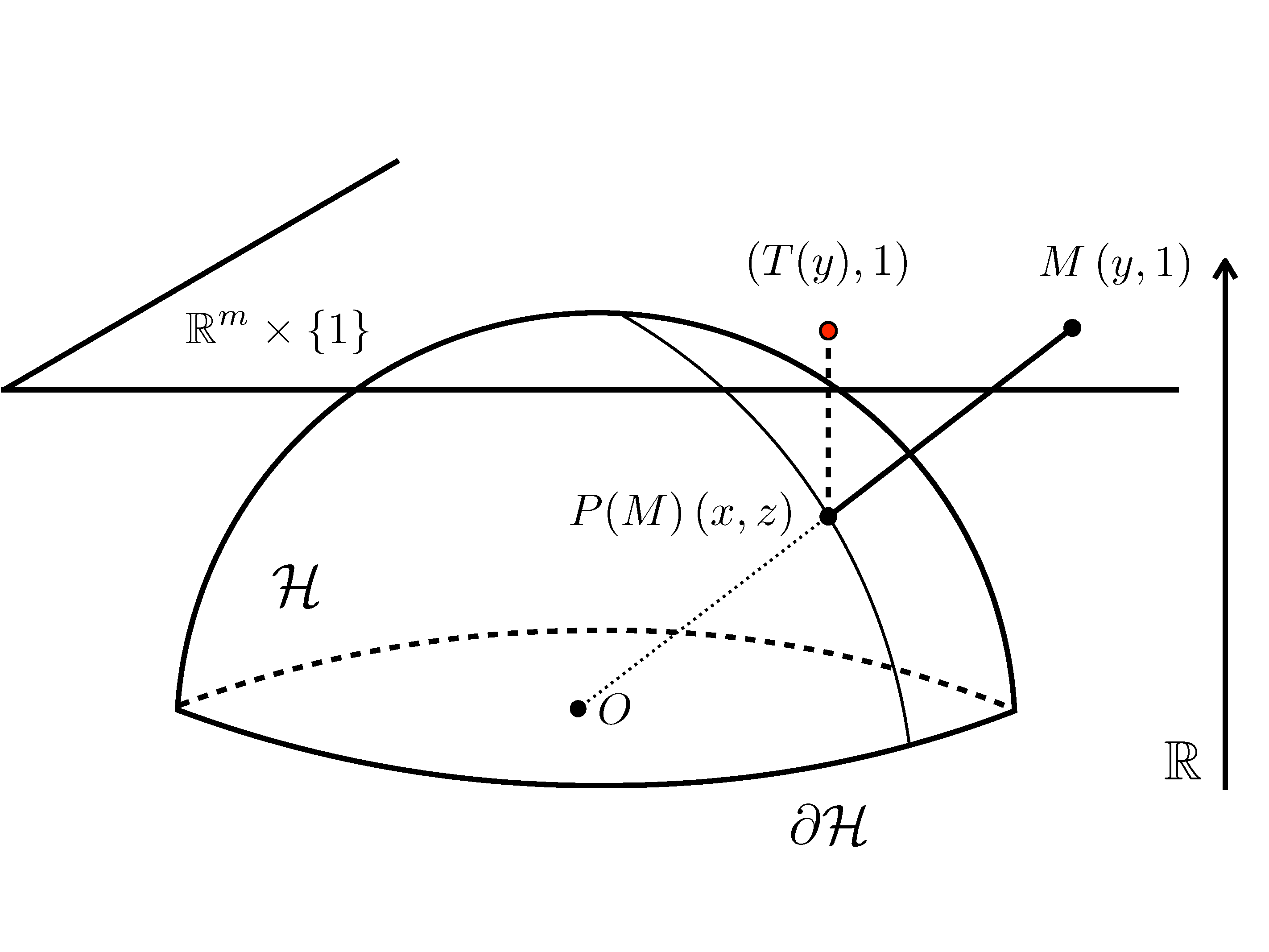}
\caption{The Poincar\'e compactification.}
\label{fig:Poincare}
\end{figure}
The Poincar\'{e} compactification is geometrically interpreted as follows (see, e.g., \cite{bib:H}).
We define the $m$-dimensional Poincar\'e hemisphere as
\[
	\mathcal{H}:=\left\{(x,z)\in \mathbb{R}^m \times\mathbb{R} \mid \|x\|^2+z^2=1,~z>0\right\}.
\]
Identifying $\mathbb{R}^m$ with $\mathbb{R}^m\times\{1\}\subset \mathbb{R}^m\times\mathbb{R}$, the space
$\mathbb{R}^m$ can be regarded as the tangent space at the north pole of $\mathcal{H}$.
A point $M=(y,1)$ on $\mathbb{R}^m\times \{1\}\cong \mathbb{R}^m$ has a one-to-one correspondence with the point $P(M)=(x,z)$ on $\mathcal{H}$ defined by 
\begin{equation*}\label{eqn:transform}
	x=\frac{y}{\sqrt{1+\|y\|^2}},\quad z=\frac{1}{\sqrt{1+\|y\|^2}},
\end{equation*}
as the intersection of $\mathcal{H}$ and the line segment in $\mathbb{R}^{m+1}$ connecting the origin in $\mathbb{R}^{m+1}$ and $M$.
The projection of $P(M)$ onto $\mathbb{R}^m\times\{1\}$ is given by $(T(y),1)$, which determines the Poincar\'{e} compactification $T_{Poin}$.
See Fig. \ref{fig:Poincare}.

\bigskip
The gradient of $\kappa$ is 
\begin{equation*}
\nabla \kappa(y) = \frac{y}{\sqrt{1+\|y\|^2}} = x,
\end{equation*}
and one can easily see that $T_{Poin}$ satisfies conditions (A0) - (A3) in Definition \ref{dfn-adm}. 
Thus, the Poincar\'{e} compactification is admissible.
Note that $T_{Poin}$ is also a {\em radial} compactification (cf. \cite{bib:Elias-Gingold}), because $\kappa$ depends only on $\|y\|$.
The normalized vector field corresponding to (\ref{normalized}) is
\begin{equation}
\label{normalized-poincare}
\frac{dx}{d\tau} = \tilde f(x, \kappa) - \langle x, \tilde f(x, \kappa)\rangle x \equiv g(x).
\end{equation}

In practical computations, we calculate the Jacobian matrix of $g$ in (\ref{normalized-poincare}) at points on $\partial \mathcal{D}$.
There are several cases where the Jacobian matrix on $\partial \mathcal{D}$ becomes discontinuous while $g$ is continuous on $\overline{\mathcal{D}}$, which are discussed with concrete examples in Section \ref{section-numerical}.
In such cases, the proposed method cannot be applied.
To avoid this situation, we may choose different compactifications, such as those described below. 
\subsection{Parabolic compactification}
\label{section-parabolic}
\begin{figure}[h]
\centering
\includegraphics[width=6.0cm]{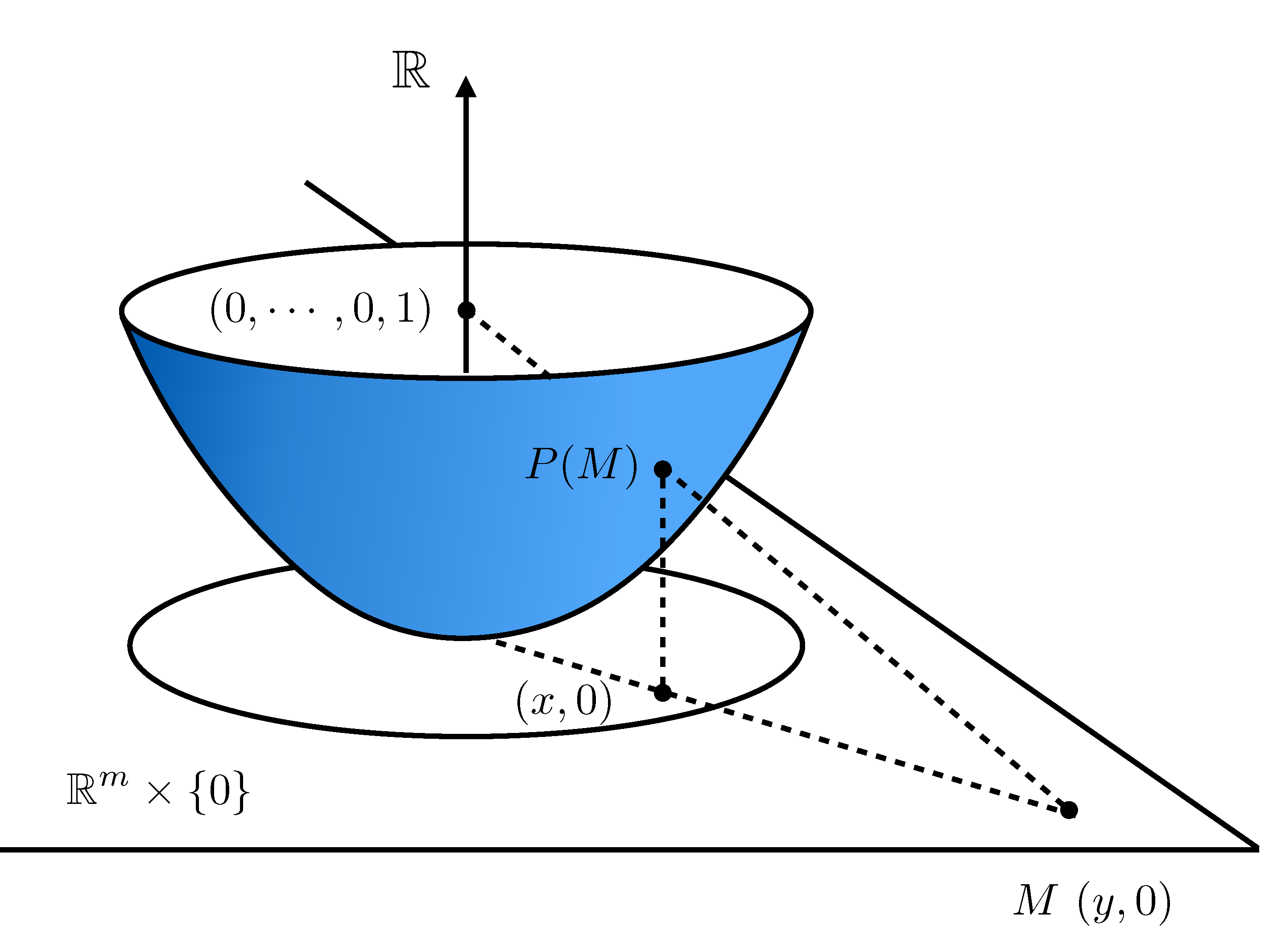}
\caption{The parabolic compactification.}
\label{fig:parabolic}
\end{figure}

A compactification using the {\em parabolic surface} $\{x_{m+1} = x_1^2 + \cdots + x_m^2\}$ is, like Poincar\'{e} compactification, an example of a radial compactification (\cite{bib:Elias-Gingold}, Section 3).
Any point $M = (y,0) = (y_1,\cdots, y_m, 0)\in \mathbb{R}^{m+1}$ is in one-to-one correspondence with the point $P(M)$, as the intersection of the parabolic surface and the line segment connecting $M$ and $(0,\cdots, 0, 1)\in \mathbb{R}^{m+1}$. 
The projection of $P(M)$ onto $\mathbb{R}^m\times \{0\}$ gives the point $(x,0) = (x_1,\cdots, x_m, 0)\in \mathbb{R}^{m+1}$, which determines the {\em parabolic compactification} $T_{para}(y) = x$.
See Fig. \ref{fig:parabolic}.
The open set $\mathcal{D} = T_{para}(\mathbb{R}^m)$ is then given by $\mathcal{D} = \{(x_1,\cdots, x_m)\in \mathbb{R}^m \mid \sum_{i=1}^m x_i^2 < 1\}$.

On the parabolic surface, the correspondence between $y\in \mathbb{R}^m$ and $x\in \mathcal{D}$ becomes
\begin{equation}
\label{parabolic}
y_j = \frac{x_j}{1-R^2},\quad R^2 = \sum_{j=1}^m x_j^2,
\end{equation}
and its inverse for $R<1$ determines $T_{para}(y) = x$, which is given concretely by 
\begin{equation}
\label{parabolic-inverse}
x_i = \frac{2y_i}{1+ \sqrt{1+4\sum_{j=1}^m y_j^2}},\quad i=1,\cdots, m.
\end{equation}
Notice that there is no square root term in (\ref{parabolic}), unlike in the Poincar\'{e} compactification.

The vector field corresponding to (\ref{transformed}) via (\ref{parabolic-inverse}) is given by
\begin{equation}
\label{normalized-parabolic-aux}
\frac{dx}{dt} = (1-R^2)^{1-d} \left[ \tilde f(x) - \frac{2\langle x,\tilde f(x)\rangle}{1+R^2}x \right],
\end{equation}
where
\begin{equation*}
\tilde f(x_1,\cdots, x_m) = (1-R^2)^d f\left(\frac{x_1}{1-R^2}, \cdots, \frac{x_m}{1-R^2}\right) = (1-R^2)^d f(y_1,\cdots, y_m).
\end{equation*}
The final equality holds from the definition (\ref{parabolic}).
The change of time variable
\begin{equation}
\label{time-transform-parabolic}
\frac{d\tau}{dt} = \frac{1}{(1-R^2(t))^{d-1} (1+R^2(t))},
\end{equation}
or
\begin{equation*}
\tau = \int_0^t \frac{ds}{(1-R^2(s))^{d-1} (1+R^2(s))},
\end{equation*}
is applied to the normalization of (\ref{normalized-parabolic-aux}).
Finally, we obtain the normalized vector field
\begin{equation}
\label{normalized-parabolic}
\frac{dx}{d\tau} = (1+R^2) \tilde f(x) - 2\langle x,\tilde f(x)\rangle x\equiv g(x).
\end{equation}
Note that the right-hand side of (\ref{normalized-parabolic}) is a polynomial, because we assumed that $f$ is a polynomial.
\section{Validating blow-up solutions}
\label{section-validation}
In this section, we provide a method for the rigorous numerical validation of blow-up solutions.
Thanks to the arguments given in Section \ref{section-compact}, we now have explicit representations of divergent solutions via compactifications of spaces and dynamical systems.
It remains to consider how we distinguish blow-up solutions from grow-up solutions for compactified systems.
This is solved by applying {\em Lyapunov functions} that are defined around the critical points at infinity.
We begin by describing our validation scenario for blow-up solutions.
Next, we review the validation method of locally defined Lyapunov functions and complete our method.
\subsection{Scenario for validating blow-up solutions}
As detailed in Section \ref{section-compact}, admissible compactifications provide a technique in dynamical systems for representing blow-up solutions.
Now, we propose the following strategy using standard methodologies of dynamical systems with admissible compactifications.
\begin{enumerate}
\item For an admissible compactification $T$, we explicitly write the normalized vector field (\ref{normalized}).
\item Find a critical point at infinity of (\ref{eqn:ODE}) in the direction $x_\ast$, i.e., a solution $x_\ast \in \partial \mathcal{D}$ of (\ref{eq-critical}). 
\item For a point $x_0\in \mathcal{D}$, validate a global trajectory $x(\tau)$ with $x(0)=x_0$ for (\ref{normalized}) that is asymptotic to $x_\ast$ as $\tau \to \infty$.
\item Compute $t_{\max}$ from (\ref{time-transform}) as 
\begin{equation}
\label{blow-up-time}
t_{\max} = \int_0^\infty \frac{d\tau}{\kappa(T^{-1}(x(\tau)))^{d-1}},
\end{equation}
and prove that $t_{\max} < \infty$.
\end{enumerate}

Step 1 can be achieved using direct calculations.
Step 2 can be performed using standard numerical methods, such as the (interval) Newton or the Krawczyk methods, or topological methods such as isolating blocks (see, e.g., \cite{bib:Mat2, bib:Mat}). 
The validity of Step 3 stems from the dynamical properties of (\ref{normalized}) on $\overline{\mathcal{D}}$.
Proposition \ref{prop-asymptotic} indicates that a solution of (\ref{eqn:ODE}) that diverges in specific direction $x_\ast$ corresponds to the solution of (\ref{normalized}) that tends to the critical point at infinity $x_\ast$. 
Because critical points at infinity correspond to equilibria for (\ref{normalized}), solutions should take an infinite time to reach them, by the standard continuation theory of differential equations.
Therefore, the statement in Step 3 is valid for describing diverge solutions from the viewpoint of dynamical systems.
In the standard theory of dynamical systems, such solutions are considered to be branches of {\em connecting orbits}.
There are many existing results for validating connecting orbits with computer assistance (see, e.g., \cite{bib:Mat2, bib:MTKO2015, bib:W}), and hence the statement in Step 3 can be reasonably executed.

However, the validation of such global solutions does not directly imply the validation of blow-up solutions.
Note that in our validations, all of the computations for solving differential equations are performed in the $\tau$-scale, which gives us no a priori information regarding solutions in $t$-scale. 
As Proposition \ref{prop-asymptotic} indicates that not just blow-up solutions but also grow-up solutions can correspond to global solutions in Step 3.
Thus, Step 4 is crucial, and represents the key point in our methodology for proving that our validated solutions are actually blow-up solutions.
The integral (\ref{blow-up-time}) contains all of the information regarding the asymptotic behavior of a given trajectory.
Thus, we must extract the asymptotic behavior of solutions in as simple a form as possible.

\begin{figure}[htbp]
\begin{minipage}{0.5\hsize}
\centering
\includegraphics[width=6.0cm]{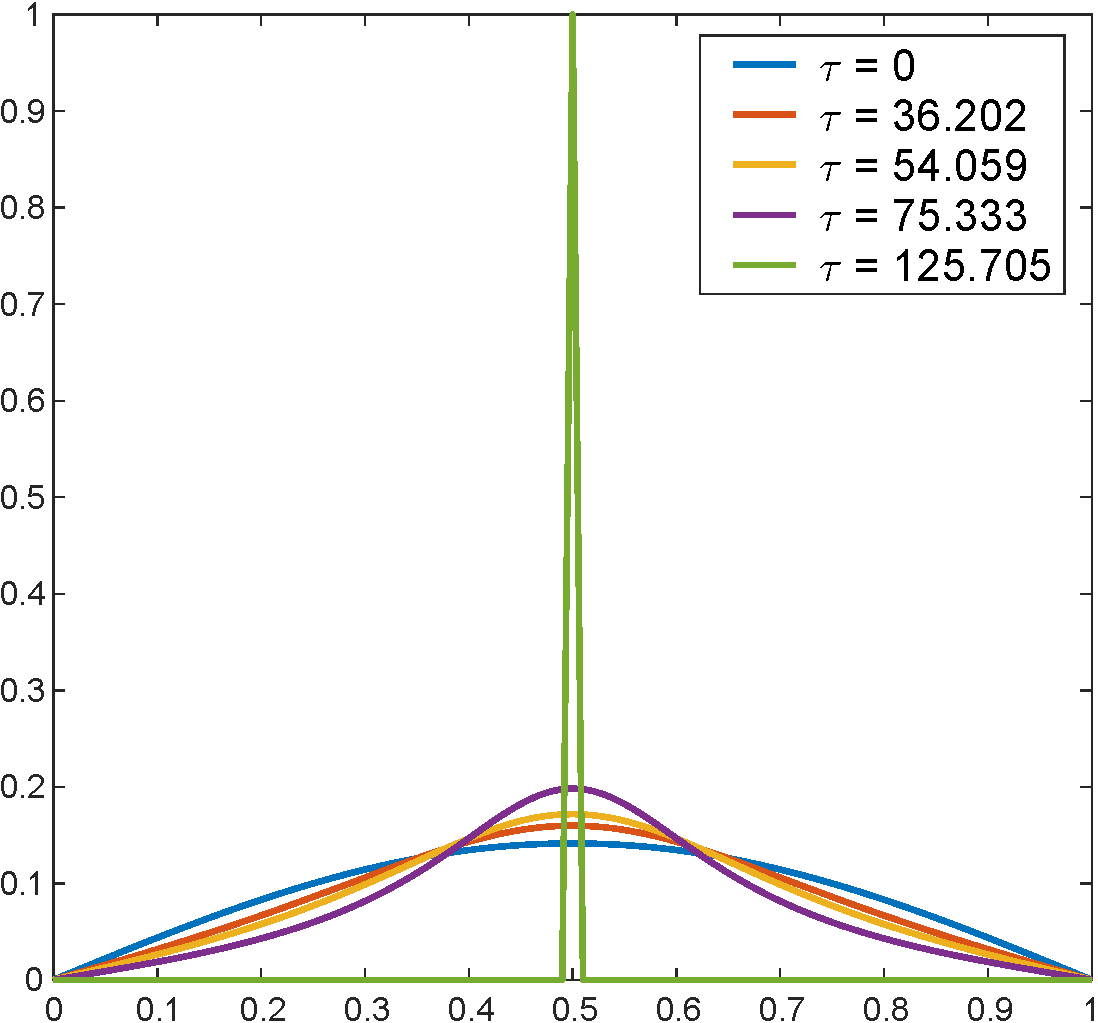}
(a)
\end{minipage}
\begin{minipage}{0.5\hsize}
\centering
\includegraphics[width=6.0cm]{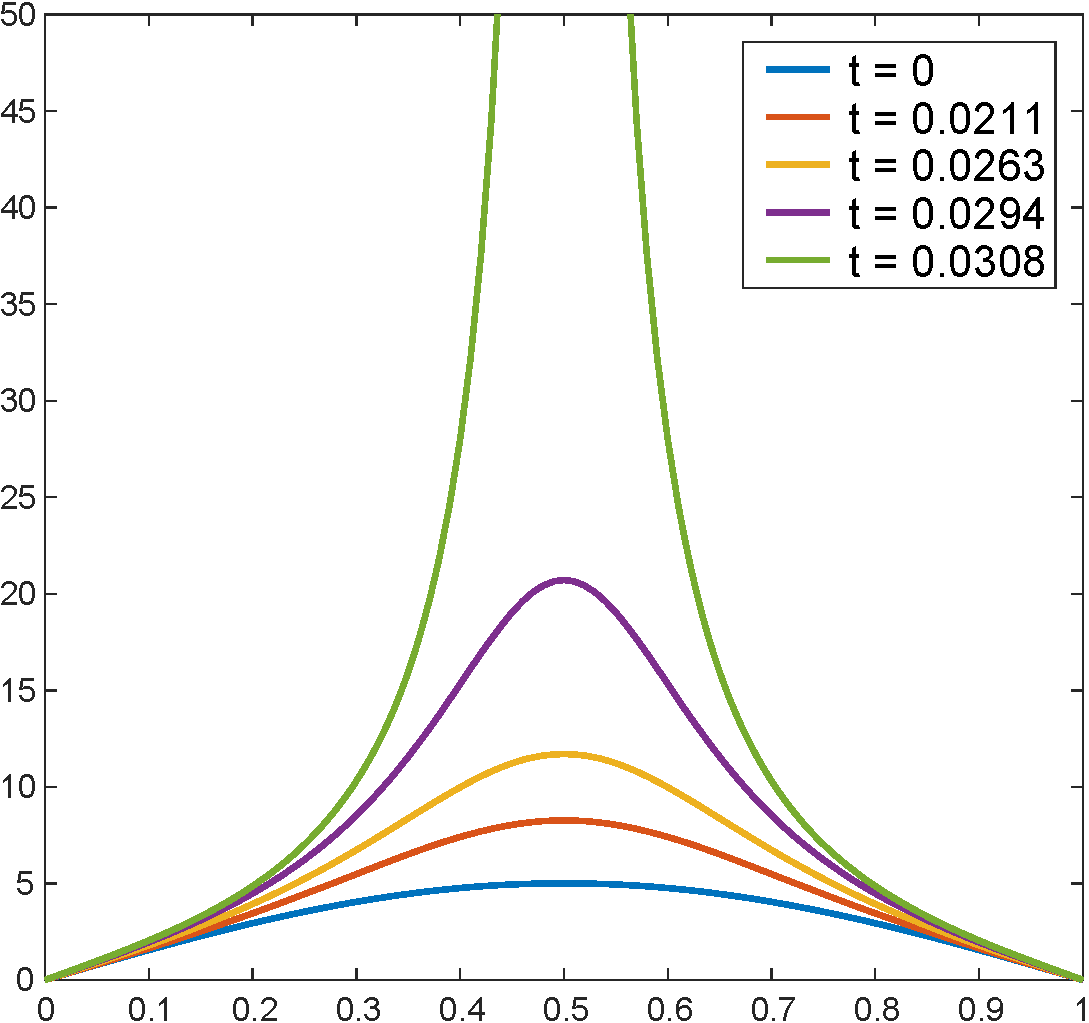}
(b)
\end{minipage}
\caption{Monotonous behavior of blow-up solutions}
\label{fig:blow-up-monotone}
(a) : Numerical solutions of (\ref{normalized}) associated with the discretized cubic nonlinear heat equation (described in Section \ref{section-heat3}) that converge to a critical point at infinity. \\
(b) : The preimage of numerical solutions in (a) in the $t$-timescale via the Poincar\'{e} compactification and the normalization (\ref{time-transform}).
The initial data is chosen as $u(x) = 5\sin (\pi x)$. 
One sees that typical ``numerical blow-up solutions'' behave asymptotically with the monotonicity of, say, the supremum norm.
\end{figure}

In many previous studies concerning the theoretical and numerical analysis of blow-up solutions \cite{bib:ALM998, bib:CHO2007, bib:N1975, bib:NU1977, bib:SS2015}, 
it is noted that typical blow-up solutions behave in a {\em monotone} way for $t$ near to blow-up times, as illustrated in Fig. \ref{fig:blow-up-monotone}.

\bigskip
We focus on such monotonicity, and provide a method of guaranteeing the monotonous behavior of solutions near to blow-ups; namely, critical points at infinity.
Our solution for this problem is the application of {\em Lyapunov functions}.

Lyapunov functions, as described in Section \ref{sec:constLyapunovfunction}, are functionals whose values monotonously decrease along solution orbits. When a dynamical system admits a Lyapunov function, it helps us to describe the global dynamics in a simple manner.
Although Lyapunov functions are usually supposed to be defined globally, the second author and his collaborators have proposed a systematic method of constructing Lyapunov functions as quadratic forms that are {\em locally defined around equilibria}. 
Discussions in \cite{bib:MHY} indicate that we can validate the explicit domain where the Lyapunov function is defined using computer assistance.
Because critical points at infinity are equilibria of transformed dynamical systems, it is natural to construct Lyapunov functions that are locally defined around critical points at infinity if we want to study the dynamics around them.
We will show that the monotonicity of Lyapunov functions describes the monotonous behavior of blow-ups well, and provides us with a new insight regarding blow-up solutions, from the viewpoint of dynamical systems.
\subsection{Lyapunov functions}
\label{sec:constLyapunovfunction}
Here, we present the constructive method for Lyapunov functions given in \cite{bib:MHY}, which enables us to describe the asymptotic behavior of solutions near to $x_\ast$. 
\begin{dfn}\rm
Let $N$ be a compact subset of $\mathbb{R}^m$, and let $U$ be an open neighborhood of $N$ in $\mathbb{R}^m$.
For a flow $\varphi$ on $\mathbb{R}^m$,
we say that a $C^1$-functional $L:U\to\mathbb{R}$ is a Lyapunov function on $N$ if $L$ satisfies the following:
\begin{enumerate}
\item 
\begin{equation*}
\left. \frac{d}{dt}L \left(\varphi(t,x)\right)\right|_{t=0}\leq 0\quad \text{ for all }x\in N;
\end{equation*}
\item If $\left.\frac{d}{dt}L\left(\varphi(t,x_\ast)\right)\right|_{t=0}=0$ holds for some $x_\ast \in N$,
then $x_\ast$ is an equilibrium of $\varphi$.
\end{enumerate}
\end{dfn}
The above properties should be usually satisfied in the whole space $\mathbb{R}^m$ or Hilbert spaces, as described in standard references regarding dynamical systems or PDEs, and hence the construction of {\em globally defined} Lyapunov functions is not an easy task, except for several typical examples.
Nevertheless, by focusing on the local dynamics around equilibria we can systematically construct a locally defined Lyapunov function of the quadratic form.
A general result for constructing a Lyapunov function around equilibria of (\ref{normalized}) is the following.
\begin{prop}[\cite{bib:MHY}, cf. \cite{bib:Mat2, bib:Mat, bib:W}]
\label{prop:Lyapunov}
Let $x_\ast$ be an equilibrium of \eqref{normalized}. 
In addition, let $N\subset\overline{\mathcal{D}}$ be a closed star-shaped subset of $\mathbb{R}^m$ with $x_\ast\in N$.
Let $Dg(x)$ be the Jacobian matrix of $g(x) = \tilde f(x) - \langle \nabla \kappa, \tilde f(x)\rangle x$ at $x$.
If $Dg$ is continuous in a neighborhood of $N$ and there exists an $m$-dimensional symmetric matrix $Y$ such that
\begin{equation}\label{eqn:ax}
	A(x):=Dg(x)^TY+YDg(x)
\end{equation}
is negative definite for all $x\in N$, then the function
\begin{equation}\label{eqn:Lyapunov-YMH}
	L(x):=(x-x_\ast)^TY(x-x_\ast)
\end{equation}
is a Lyapunov function on $N$.
We shall call the domain $N$ a {\em Lyapunov domain} (for $L$).
\end{prop}
The negative definiteness of the matrix $A(x)$ implies that all eigenvalues of $A(x)$ have nonzero real parts. 
As we will see below, the matrix $Y$ contains the information regarding the diagonalization of $Dg(x_\ast)$.
The explicit expression of $L$ in (\ref{eqn:Lyapunov-YMH}) and its properties imply that $x_\ast$ is a {\em unique} equilibrium in $N$.

\bigskip
Following \cite{bib:MHY}, we can explicitly construct $Y$ as follows.
Let $\Lambda=V^{-1}Dg(x_\ast)V$ be an (approximate) diagonalization of $Dg(x_\ast)$.
Assume that no eigenvalues $\{\lambda_j\}_{j=1}^{m}$ of $\Lambda$ lie on the imaginary axis. 
Then, we define $I^*=\mathrm{diag}(i_1,...,i_m)$, $\hat Y=V^{-H}I^*V^{-1}$, and $Y=Re(\hat Y)$,
where
\[
	i_j=\left\{
	\begin{array}{rl}
	1&(Re(\lambda_j)<0),\\
	-1&(Re(\lambda_j)>0)
	\end{array}
	\right.
\]
and $V^{-H}$ is the Hermitian transpose of $V^{-1}$.

In particular, the above procedure becomes slightly simpler for {\em stable} critical points at infinity, which encompasses the cases we treat below. 
Indeed, let $x_\ast$ be a critical point at infinity ($g(x_\ast)=0$), and assume that it holds that $Re(\lambda_j) < 0$ for all $j=1,\cdots, m$.
Then, the matrix $I^\ast$ is defined as the identity matrix, and hence $\hat Y=V^{-H}V^{-1}$.
Throughout the remainder of our arguments, this property regarding eigenvalues central to our considerations.
\begin{rem}
Lyapunov functions themselves can be applied not only to equilibria that are stable, but also to those of saddle-type.
It has been observed that the condition stated in Proposition {\rm\ref{prop:Lyapunov}} is equivalent to a sufficient condition of {\em cone conditions} (\cite{bib:ZCov}. See also \cite{bib:Mat2, bib:MHY}). In other words, the stable and unstable manifolds of equilibria can be given by graphs of smooth functions. We omit the details, because we only deal with the simplest case, that is, the case that $N$ is a subset of the stable manifold of $x_\ast$.
\end{rem}
\subsection{Enclosing the blow-up time}
We are now ready to present our method for enclosing the blow-up time $t_{\max}$ given by (\ref{blow-up-time}).
Let $x_\ast\in\overline{\mathcal{D}}$ be a critical point at infinity.
Using the technique described in the previous subsection, we can validate a Lyapunov function $L$ on a given compact domain $\tilde N\subset \overline{\mathcal{D}}$ containing $x_\ast$.
Recall that a typical property of Lyapunov functions is that they can validate the monotonous behavior of dynamical systems.
We can apply this property to the {\em re-parameterization} of dynamical systems around equilibria.

Let $x=x(\tau)$ be a solution orbit for the normalized vector field (\ref{normalized}), with $x(0) = x_0 \in \tilde N$. 
For a given solution orbit with fixed initial data $x_0$, the time $\tau$ and the value of $L$ correspond one-to-one as long as $x(\tau)\in \tilde N$, 
because $L$ monotonously decreases along solution orbits.
Of course, $L(x(\tau))$ is smooth in $\tau$, because $x$ smoothly depends on $\tau$.
This property tells us that we can apply a technique called {\em Lyapunov tracing}, as proposed in \cite{bib:MHY}, which is a smooth change of variables between $\tau$ and $L$ along $x(\tau)$ such that the integral in (\ref{blow-up-time}) can be estimated by an $L$-integral.
If the solution orbit $x(\tau)$ is contained in $\tilde N$ for all $\tau\geq 0$, 
then the change of coordinates allows us to estimate the integral (\ref{blow-up-time}) explicitly.

\begin{figure}[h]
{\centering
\includegraphics[width=6.0cm]{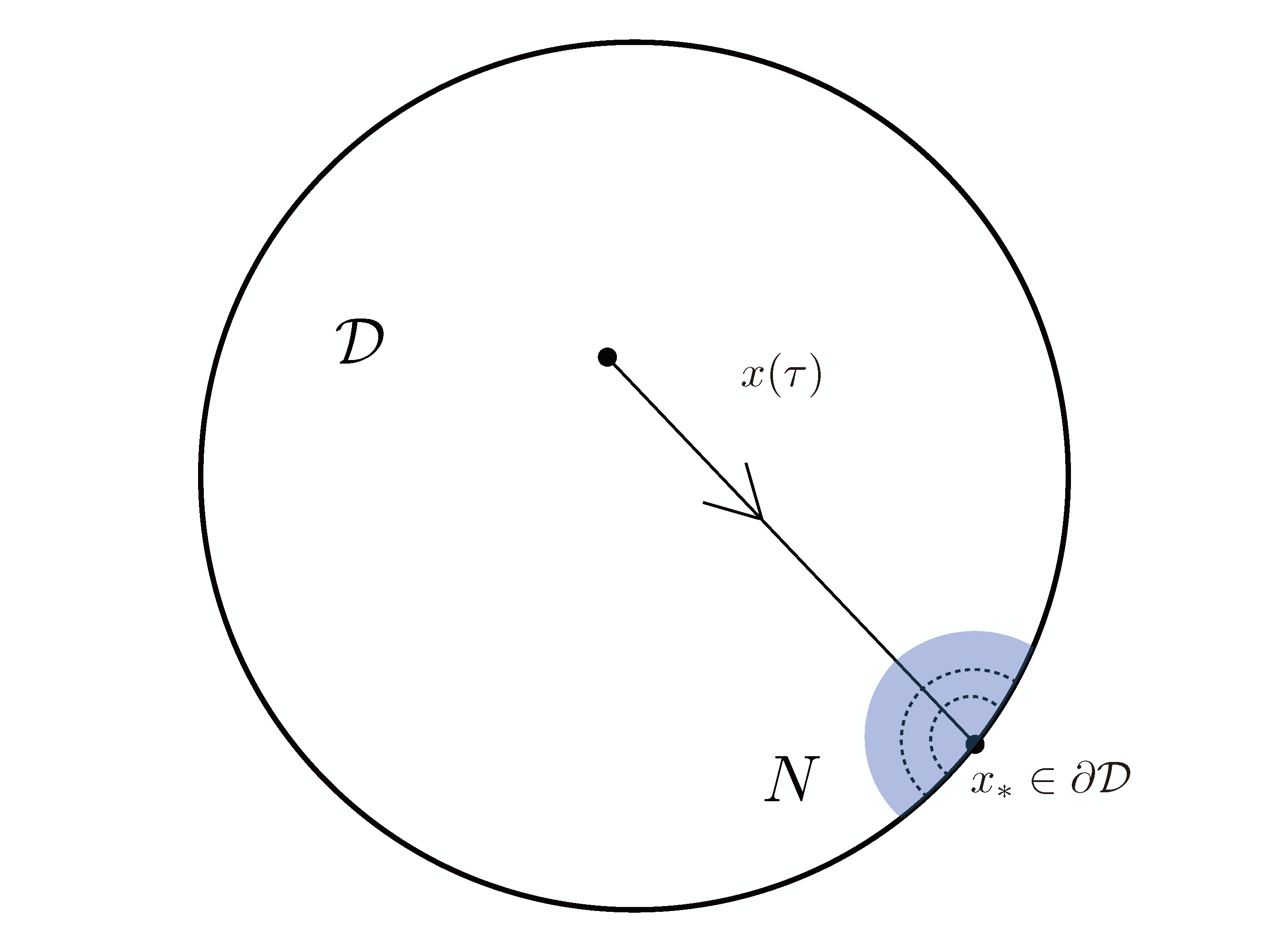}
\caption{Behavior of a solution $x(\tau)$.}
\label{fig:behavior}}
Dotted curves represent level sets of the Lyapunov function $L$ of a stable critical point at infinity $x_\ast$.
Once trajectories enter $N$, they remain at $N$ under the forward time flow and tend to $x_\ast$, so that they intersect each level set of $L$ transversely.
\end{figure}

\bigskip
With this observation in mind, we present a concrete procedure for estimating $t_{\max}$.
Note that the essential ideas applied for { our} estimate below also apply for other compactifications.
Let $L$ be a Lyapunov function of the form (\ref{eqn:Lyapunov-YMH}) on an appropriate set $\tilde N\subset \overline{\mathcal{D}}$ that contains $x_\ast$\footnote{The validation of a Lyapunov domain on the basis of interval (or affine) arithmetic is performed on a rectangular domain $\mathcal{N}\subset \mathbb{R}^m$ around $x_\ast$. 
This usually contains a region in $\mathbb{R}^m\setminus \overline{\mathcal{D}}$, 
while the set being a Lyapunov domain makes sense only in $\mathcal{N}\cap \overline{\mathcal{D}}$,
because the dynamics \eqref{normalized} is valid only in $\overline{\mathcal{D}}$.}.
Assume that there exists a constant $\delta_0 >0$ such that $Re(\lambda_j)\leq -\delta_0 < 0$ for all eigenvalues $\lambda_j(x)$ of $Dg(x)$, with $x\in \tilde N$. 
Choose a subset $N\subset\overline{\mathcal{D}}$ of the form
\begin{equation}\label{eqn:domainofLyapunov}
	N=\left\{x\in\overline{\mathcal{D}}:L(x)\le\varepsilon^2\right\},
\end{equation}
so that $N\subset \tilde N$, where $\varepsilon>0$.
The set $N$ is then an attracting neighborhood of $x_\ast$ in $\overline{\mathcal{D}}$. That is, it holds that $\langle g(x), \nu(x) \rangle<0$ for all $x\in \overline{\partial N\setminus \partial \mathcal{D}}$, where $\nu(x)$ is an outer unit normal vector to $\partial N$ at $x\in \partial N$ (cf. \cite{bib:Mat}).
This directly follows from the definitions of $L$ and $N$.
Moreover, there exist $c_{\tilde N},~c_1>0$ such that
\begin{align}\label{eqn:Lyapieq}
	\frac{dL}{d\tau}(x(\tau))&<-c_{\tilde N}\|x-x_\ast\|^2\nonumber\\
	&\le-c_{\tilde N}c_1L(x).
\end{align}
\begin{rem}
The constant $c_1>0$ is given by
\[
	c_1 = \max_{j=1,...,m}\{\mu_j^{-1}:\mu_j~\mbox{is an eigenvalue of}~Y\}.
\]
It can be shown that all eigenvalues of $Y$ are positive in our case (Proposition {\rm 2} in {\cite{bib:MHY}}).
For example, $\tilde N$ can be chosen as
\[
	\left\{x\in \overline{\mathcal{D}}:|x_i-x_{\ast,i}|\le\varepsilon \sqrt{c_1},~i=1,...,m\right\}.
\]
The constant $c_{\tilde N}$ is then given by
\[
	c_{\tilde N}=\min_{\substack{x\in \tilde N\\j=1,...,m}}\{|\lambda_j(x)|:\lambda_j(x)~\mbox{is an eigenvalue of}~A(x)\}.
\]
Note that $c_1$ and $c_{\tilde N}$ can be rigorously estimated using validated computations and, say, the Gerthgorin Circle Theorem.
\end{rem}
If we have the solution $x:[0,\tau_N]\to\mathcal{D}$ of \eqref{normalized} satisfying $x(\tau_N)\in N$,
then $x$ corresponds to a solution $y:[0,t_N]\to \mathbb{R}^m$ of \eqref{eqn:ODE},
where $t_N$ is given by
\[
	t_N:=\int_0^{\tau_N}\frac{d\tau}{\kappa\left(T^{-1}\left(x(\tau)\right)\right)^{d-1}}.
\]
Furthermore,
the solution of \eqref{normalized} with the initial point $x(\tau_N)$ asymptotically goes to $x_\ast$ as $\tau\to\infty$,
so that it decreases the Lyapunov function $L$ monotonically, as illustrated in Fig. \ref{fig:behavior}.
We also note that the property $\{x(\tau)\}_{\tau\in [\tau_N,\infty)}\subset N$ is guaranteed to hold, because $N$ is an attracting neighborhood of $x_\ast$.

The solution $y$ of \eqref{eqn:ODE} corresponding to $x:[\tau_N,\infty)\to\overline{\mathcal{D}}$ has the interval of existence $[t_N,t_{\max})$,
where
\begin{align*}
	t_{\max}
	&=t_N+\int_{\tau_N}^\infty\frac{d\tau}{\kappa\left(T^{-1}\left(x(\tau)\right)\right)^{d-1}}.
\end{align*}
Because $x(\tau_N)\in\mathcal{D}$, $t_N$ is finite.
Thus, we can prove that $t_{\max}<\infty$ (i.e., $y$ is a blow-up solution) if
\[
	t_{\max}-t_N=\int_{\tau_N}^\infty\frac{d\tau}{\kappa\left(T^{-1}\left(x(\tau)\right)\right)^{d-1}}
\]
is finite.
Explicit estimates of $t_{\max}-t_N$ depend on the choice of admissible compactifications. 
We provide explicit estimates in the case of the Poincar\'{e} and the parabolic compactifications.
\subsubsection{Case 1: Poincar\'{e} compactifications} 
First, we provide an explicit estimate of $t_{\max}-t_N$ under the Poincar\'{e} compactification. 
From the fact that $x=T(y)=y/\kappa(y)$, it holds that
\[
	\|y\|^2=\frac{\|x\|^2}{1-\|x\|^2},
\]
which yields that
\[
	\kappa\left(T^{-1}(x)\right)=\left(1+\frac{\|x\|^2}{1-\|x\|^2}\right)^{\frac{1}{2}}=\left(1-\|x\|^2\right)^{-\frac{1}{2}}.
\]
Because $\|x_\ast\|^2=1$, we have that
\begin{align*}
	1-\|x\|^2&=2\left(\langle x_\ast, x_\ast-x\rangle\right)-\|x-x_\ast\|^2
	\le2\|x-x_\ast\|\le 2(c_1L)^{\frac{1}{2}}.
\end{align*}
This inequality yields that
\begin{align}
	t_{\max} - t_N 
	&= \int_{\tau_N}^\infty (1-\|x(\tau)\|^2)^{\frac{d-1}{2}}d\tau\nonumber\\
	&\leq \int_{\tau_N}^\infty (2\|x(\tau)-x_\ast\|)^{\frac{d-1}{2}}d\tau\nonumber\\
	&\leq \int_{\tau_N}^\infty (2(c_1L(x(\tau)))^{1/2})^{\frac{d-1}{2}}d\tau\nonumber\\
	&\leq \frac{2^{\frac{d-1}{2}}}{c_{\tilde N}} \int_0^{L(x(\tau_N))} (c_1L)^{\frac{d-1}{4}-1}dL\nonumber\\
	&= \frac{2^{\frac{d-1}{2}} c_1^{\frac{d-5}{4}} }{c_{\tilde N}} \frac{4}{d-1}L(x(\tau_N))^{\frac{d-1}{4}},\label{eqn:Blowuptime}
\end{align}
where $c_1$ and $c_{\tilde N}$ are constants in \eqref{eqn:Lyapieq}.
The re-parameterization in the above estimates is realized along the solution $\left\{x(\tau)\right\}_{\tau \in [\tau_N,\infty)}$.
Note that it holds that $x(\tau)\in N$ for any $\tau\in[\tau_N,\infty)$, because $N$ is an attracting neighborhood of $x_\ast$.
Furthermore, the value $L\left(x(\tau_N)\right)$ is explicitly evaluated by
\[
	L\left(x(\tau_N)\right)=\left(x(\tau_N)-x_\ast \right)^TY\left(x(\tau_N)- x_\ast\right).
\]
\begin{rem}
We use the numerical verification library called {\em kv} \cite{bib:kashikashi,bib:kv} to calculate an enclosure of the solution $x:[0,\tau_N]\to \mathcal{D}$ of \eqref{normalized-poincare} that satifies $x(\tau_N)\in N$ as well as $t_N$, where
\begin{align}
	t_N
	=\int_0^{\tau_N}\frac{d\tau}{\kappa\left(T^{-1}\left(x(\tau)\right)\right)^{d-1}}
\label{tN_poincare}
	=\int_0^{\tau_N}\left(1-\|x(\tau)\|^2\right)^{\frac{d-1}{2}}d\tau.
\end{align}
\end{rem}
Summarizing the above arguments, we obtain rigorous upper and lower bounds of blow-up times.
\begin{thm}\label{th:main_poincare}
Let $x_\ast\in \partial \mathcal{D}$ be an equilibrium of \eqref{normalized-poincare}.
Assume that there is a compact neighborhood $\tilde N$ of $x_\ast$ in $\overline{\mathcal{D}}$ that admits a Lyapunov function $L$ of the form \eqref{eqn:Lyapunov-YMH} with $x_\ast$, and that for all eigenvalues of $Dg(x)$ with $x\in \tilde N$, $Re(\lambda_j(x)) \leq -\delta_0 < 0$ holds for some $\delta_0 > 0$.
Let $N = \{x\in \overline{\mathcal{D}} : L(x)\leq \varepsilon^2\}$, where $\varepsilon > 0$ is chosen such that $N \subset \tilde N$.
Furthermore, let $\{x(\tau)\}_{\tau\in [0,\tau_N]}$ be a solution orbit of \eqref{normalized-poincare} with $x(0) = x_0$, such that $x(\tau_N)\in {\rm int}\,N$.
Then, the mapping $y(t)$ defined by $y(t) = T^{-1}(x(\tau))$, with $\tau \in [0,\infty)$, is a solution of \eqref{eqn:ODE} with $y(0) = T^{-1}(x_0)$, whose maximal existence time $t_{\max}$ is contained in an interval $[\underline{t_{\max}}, \overline{t_{\max}}]$. Here, $\underline{t_{\max}}$ is given by $t_N$ in \eqref{tN_poincare}, and
\begin{equation*}
\overline{t_{\max}} = \underline{t_{\max}} + \frac{2^{\frac{d-1}{2}} c_1^{\frac{d-5}{4}} }{c_{\tilde N}} \frac{4}{d-1} \varepsilon^{\frac{d-1}{2}},
\end{equation*}
where the constants $c_1$ and $\tilde c_N$ are given in \eqref{eqn:Lyapieq}.
In particular, if $\overline{t_{\max}} < \infty$, then $y(t)$ is a blow-up solution with the blow-up time
$t_{\max} \in [ \underline{t_{\max}}, \overline{t_{\max}} ]$.
\end{thm}
\subsubsection{Case 2: parabolic compactifications} 
\label{section-blowup-parabolic}
Next, we provide an explicit estimate of $t_{\max}-t_N$ under parabolic compactification.

Let $x(\tau)$, with $x(0)=x_0\in \mathcal{D}$, be a global solution of (\ref{normalized-parabolic}) that is asymptotic to a critical point at infinity $x_\ast \in \partial \mathcal{D}$.
From (\ref{time-transform-parabolic}), the maximal existence time of $x(\tau)$ in $t$-time scale is then
\begin{equation}
t_{\max} = \int_0^\infty (1-R(\tau)^2)^{d-1}(1+R(\tau)^2) d\tau
\end{equation}
where $R(\tau)^2 = \sum_{i=1}^m x_j(\tau)^2$.

Assume that there is a neighborhood $\tilde N$ of $x_\ast$ in $\overline{\mathcal{D}}$ such that all eigenvalues of the Jacobian matrix $Dg$ of the vector field $g$ given in (\ref{normalized-parabolic}) have negative real parts for all $x\in \tilde N$.
We also assume that $\tilde N$ admits a Lyapunov function $L$, and contains $N \equiv \{x\in \overline{\mathcal{D}}\mid L(x)\leq \varepsilon^2\}$ for some $\varepsilon > 0$.
As in the case of the Poincar\'{e} compactification, the inequality $1-R^2 \leq 2(c_1 L)^{1/2}$ holds. 
Obviously, it also holds that $1+R^2 \leq 2$, and hence we obtain the following upper bound on the maximal existence time of $x(t)$ in the $t$-time scale: 
\begin{align}
t_{\max} - t_N &= \int_{\tau_N}^\infty (1-R(\tau)^2)^{d-1}(1+R(\tau)^2) d\tau\nonumber\\
	&\leq \int_{\tau_N}^\infty 2 (2(c_1 L)^{1/2})^{d-1}d\tau = 2^d\int_{\tau_N}^\infty  (c_1 L)^{\frac{d-1}{2}}d\tau\nonumber\\
	&\leq \frac{2^d}{c_{\tilde N}}\int_0^{L(x(\tau_N))}  (c_1 L)^{\frac{d-1}{2}-1}dL\nonumber\\
	&= \frac{2^d c_1^{\frac{d-3}{2}}}{c_{\tilde N}} \frac{2}{d-1}L(x(\tau_N))^{\frac{d-1}{2}},\label{eqn:Blowuptime_para}
\end{align}
where $x(\tau_N)$ is a point on the trajectory $x(\tau)$ at the time $\tau = \tau_N$ such that $x(\tau_N)\in N$ and
\begin{equation*}
t_N = \int_0^{\tau_N} (1-R(\tau)^2)^{d-1}(1+R(\tau)^2) d\tau.
\end{equation*}
The re-parameterization in the above estimates is realized along the solution $\{x(\tau)\}_{\tau \in [\tau_N,\infty)}$.

Summarizing these arguments, we obtain rigorous upper and lower bounds of blow-up times.
\begin{thm}\label{th:main_para}
Let $x_\ast\in \partial \mathcal{D}$ be an equilibrium of \eqref{normalized-parabolic}.
Assume that there is a compact neighborhood $\tilde N$ of $x_\ast$ in $\overline{\mathcal{D}}$ that admits a Lyapunov function $L$ of the form \eqref{eqn:Lyapunov-YMH} with $x_\ast$, and that it holds that for all eigenvalues of $Dg(x)$ with $x\in \tilde N$, $Re(\lambda_j(x))\leq -\delta_0 < 0$ holds for some $\delta_0 > 0$.
Let $N = \{x\in \overline{\mathcal{D}} : L(x)\leq \varepsilon^2\}$, where $\varepsilon > 0$ is chosen such that $N \subset \tilde N$.
Furthermore, let $\{x(\tau)\}_{\tau\in [0,\tau_N]}$ be a solution orbit of \eqref{normalized-parabolic} with $x(0) = x_0$, such that $x(\tau_N)\in {\rm int}\,N$.
Then, the mapping $y(t)$ defined by $y(t) = T^{-1}(x(\tau))$, with $\tau \in [0,\infty)$, is a solution of \eqref{eqn:ODE} with $y(0) = T^{-1}(x_0)$, whose maximal existence time $t_{\max}$ is contained in an interval $[\underline{t_{\max}}, \overline{t_{\max}}]$. Here, $\underline{t_{\max}}$ is given by $t_N$ in \eqref{tN_poincare}, and
\begin{equation*}
\overline{t_{\max}} = \underline{t_{\max}} + \frac{2^d c_1^{\frac{d-3}{2}}}{c_{\tilde N}} \frac{2}{d-1}\varepsilon^{d-1},
\end{equation*}
where the constants $c_1$ and $\tilde c_N$ are given in \eqref{eqn:Lyapieq}.
In particular, if $\overline{t_{\max}} < \infty$, then $y(t)$ is a blow-up solution with the blow-up time
$t_{\max} \in [ \underline{t_{\max}}, \overline{t_{\max}} ]$.
\end{thm}
\subsubsection{Remarks on the topology of ``blow-up'' and ``blow-up time''}
To conclude this section, we remark on the topology of blow-up phenomena and blow-up times.
Throughout our arguments, we have employed the standard Euclidean norm $\|y\| = \sqrt{y_1^2 + \cdots + y_m^2}$.
Consequently, the blow-up phenomena in our validations themselves make sense for {\em every} norm in $\mathbb{R}^m$, because all arbitrary norms in $\mathbb{R}^m$ are equivalent.
In particular, for any norm $\|\cdot \|'$ in $\mathbb{R}^m$, there exist positive constants $c$ and $C$ such that
\begin{equation*}
c\|x\| \leq \|x\|' \leq C\|x\|,\quad \forall x\in \mathbb{R}^m.
\end{equation*}

In contrast, the enclosure of blow-up times in our validations makes sense only for the Euclidean norm $\|\cdot \|$.
Indeed, the above estimates are based on the Euclidean norm.
If different norms are chosen, then in general the estimates of the validated blow-up times will change.

{
\subsection{Remarks on the validation}
In our validation method, the success of \eqref{eqn:Lyapieq} implies that the validated solution is actually a blow-up solution.
In other words, \eqref{eqn:Lyapieq} gives a sufficient condition for validating blow-up solutions, not a necessary condition.
Additionally, one of the key points is an appropriate choice of normalizations so that
\[
	t_{\max} < \infty \text{ in the original time scale},\quad \tau_{\max} = \infty\text{ in the normalized time scale}
\]
for blow-up solutions.
On the other hand, when an ODE has a grow-up solution, the maximal existence time $t_{\max}$ must be  infinity, which implies that the integral in \eqref{blow-up-time} would diverge for any normalization.


The upper bound $\overline{t_{\max}}$ in Theorem \ref{th:main_poincare} and Theorem \ref{th:main_para} is always bounded  as long as $\underline{t_{\max}}, c_1, c_{\tilde N}$ can be validated to be finite.
However, the boundedness of $\overline{t_{\max}}$ depends on the hyperbolicity of the critical point at infinity $x_\ast$.
For example, if $x_\ast$ is non-hyperbolic, namely, one of the eigenvalues of $Dg(x_\ast)$ has a zero real part, then $c_{\tilde N}$ becomes infinite.
Moreover, $\overline{t_{\max}}$ also becomes infinite if $d=1$, namely, the vector field is linear.
We also note that 
it is not obvious even for a blow-up solution whether $\underline{t_{\max}}, c_1, c_{\tilde N}$ can be computed to be finite or not in actual computations.


Our validation method in fact can be applied for any vector field $f$ in \eqref{eqn:ODE} as long as the
normalized vector field $g$ in \eqref{normalized} is $C^1$ on $\overline{\mathcal{D}}$, namely, the function $g$ on $\mathcal{D}$ can be extended $C^1$-smoothly to  $\overline{\mathcal{D}}$.
In particular, the vector field $f$ does not restricted to be polynomial.

}
\section{Numerical validation examples}
\label{section-numerical}
In this section, we present several numerical validation examples to demonstrate the applicability of our method.
All computations were carried out on Cent OS 6.3, Intel(R) Xeon(R) CPU E5-2687W@3.10 GHz using the kv library \cite{bib:kv} to rigorously compute the trajectories of ODEs.
\subsection{$dy/dt = y^2$}
As a benchmark test, we first considered the following initial value problem: 
\begin{equation}\label{eqn:ex1}
	\frac{dy}{dt}=y^2,~y(0)=a>0.
\end{equation}
The exact solution of this problem is
$y(t)=(a^{-1}-t)^{-1}$, which blows up at $t_{\max}=a^{-1}$.
We chose the Poincar\'e compactification as an admissible compactification.
The normalized problem of \eqref{eqn:ex1} corresponding to \eqref{normalized-poincare} is given by
\begin{equation}\label{eqn:ex1_normed}
	\frac{dx}{d\tau}=x^2-x^4,~x(0)=\frac{a}{\sqrt{1+a^2}}.
\end{equation}
The critical points at infinity are $x_\ast=\pm 1$.
We set $N$ as in \eqref{eqn:domainofLyapunov}, with $x_\ast=1$ and $\varepsilon = 1.0\times 10^{-10}$.
The Lyapunov function can be then constructed on the set $\tilde{N}$ containing $N$.
We rigorously computed a solution of \eqref{eqn:ex1_normed} with $a=1/4$, satisfying $x\to 1$.
When $\tau_N=15$, we enclosed
\[
	x(\tau_N)\in [0.999999999941121,0.999999999941125]\subset N
\]
and $t_N\in [3.9999713430726937,4.0000069540703374]$,
where $[\cdot,\cdot]$ denotes a real interval.
Therefore, from Theorem \ref{th:main_poincare}, we have proved that the solution $x$ with $a=1/4$ asymptotically goes to the critical point at infinity $x_\ast=1$ as $\tau\to\infty$.
From \eqref{eqn:Blowuptime}, we have that
$t_{\max}-t_N\le 1.0851585850519801\times 10^{-5}$.
Consequently, the blow-up time of the solution of \eqref{eqn:ex1} was rigorously enclosed as
\[
	t_{\max}\in[3.9999713430726937,4.0000178056561886].
\]
The execution time was about 0.104 seconds.
Indeed, this interval includes the exact blow-up time  $t_{\max}=4$.
\subsection{A two-dimensional system}
\begin{figure}[h]
{\centering
\includegraphics[width=8.0cm]{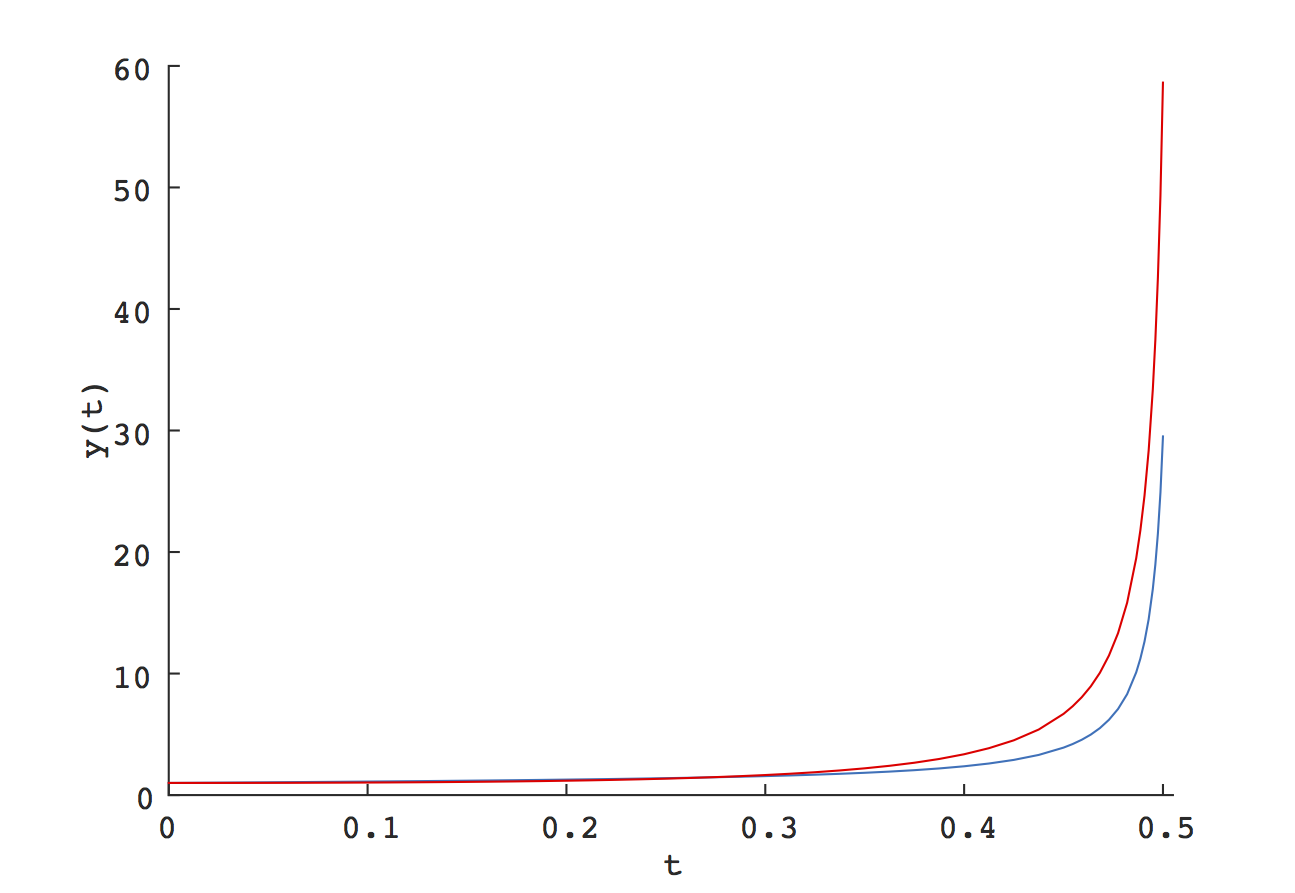}
\caption{Behavior of a solution $y(t)$ of \eqref{eqn:ex2}.}
\label{fig:behaviorex2}}
The blue curve represents $y_1(t)$, and the red curve represents $y_2(t)$.
\end{figure}

The second example is the following ODE in $\mathbb{R}^2$:
\begin{equation}\label{eqn:ex2}
\left\{
\begin{array}{l}
	\displaystyle\frac{dy_1}{dt}=y_1^2+y_2^2-1,\\[4mm]
	\displaystyle\frac{dy_2}{dt}=5(y_1y_2-1),
\end{array}
\right.
\end{equation}
with the initial value $(y_1(0),\,y_2(0))=(1,1)$.
We also chose the Poincar\'e compactification as an admissible compactification.
The normalized problem of \eqref{eqn:ex2} corresponding to \eqref{normalized-poincare} is given by
\begin{equation}\label{eqn:ex2_normed}
\left\{
\begin{array}{rl}
	\displaystyle\frac{dx_1}{d\tau}=&\hspace*{-3mm}(1-x_1^2)(2(x_1^2+x_2^2)-1) -5x_1x_2(x_1^2+x_1x_2+x_2^2-1),\\[2mm]
	\displaystyle\frac{dx_2}{d\tau}=&\hspace*{-3mm}5(1-x_2^2)(x_1^2+x_1x_2+x_2^2-1) -x_1x_2(2(x_1^2+x_2^2)-1).
\end{array}
\right.
\end{equation}
There are six critical points at infinity: $(\pm1,0)$, $(\pm1/\sqrt{5},$ and $\pm2/\sqrt{5})$.
We set $N$ as in \eqref{eqn:domainofLyapunov}, with $x_\ast=(1/\sqrt{5},2/\sqrt{5})$ and $\varepsilon = 1.0\times 10^{-10}$.
The Lyapunov function can be then constructed on the set $\tilde{N}$ containing $N$.
We rigorously computed a solution orbit of \eqref{eqn:ex2_normed} with $x(0)=(1/\sqrt{3},1/\sqrt{3})$ satisfying $x\to x_\ast$.
When $\tau_N=7$, we obtained
\[
{\small x(\tau_N)\in
\left(
\begin{array}{l}
\left[0.447213595573401, 0.447213595573404\right]\\
\left[0.894427190963148, 0.894427190963151\right]
\end{array}
\right)\subset N}
\]
and $t_N\in [0.50680733588232473,0.50681331159709609]$.
Therefore, from Theorem \ref{th:main_para}, we have proved that the solution $x$ with the initial value $x(0)$ asymptotically goes to the critical point at infinity $x_\ast$ as $\tau\to\infty$.
From \eqref{eqn:Blowuptime}, we have that $t_{\max}-t_N\le7.6274327476907032\times 10^{-6}$.
Consequently, the blow-up time of the solution of \eqref{eqn:ex2} was rigorously enclosed as
\[
	t_{\max}\in[0.50680733588232473,0.50682093902984382].
\]
The execution time was about 1.217 seconds.
The profile of this blow-up solution is given in Fig. \ref{fig:behaviorex2}.
\subsection{Spiral blow-up}
\begin{figure}[h]
{\centering
\includegraphics[width=8.0cm]{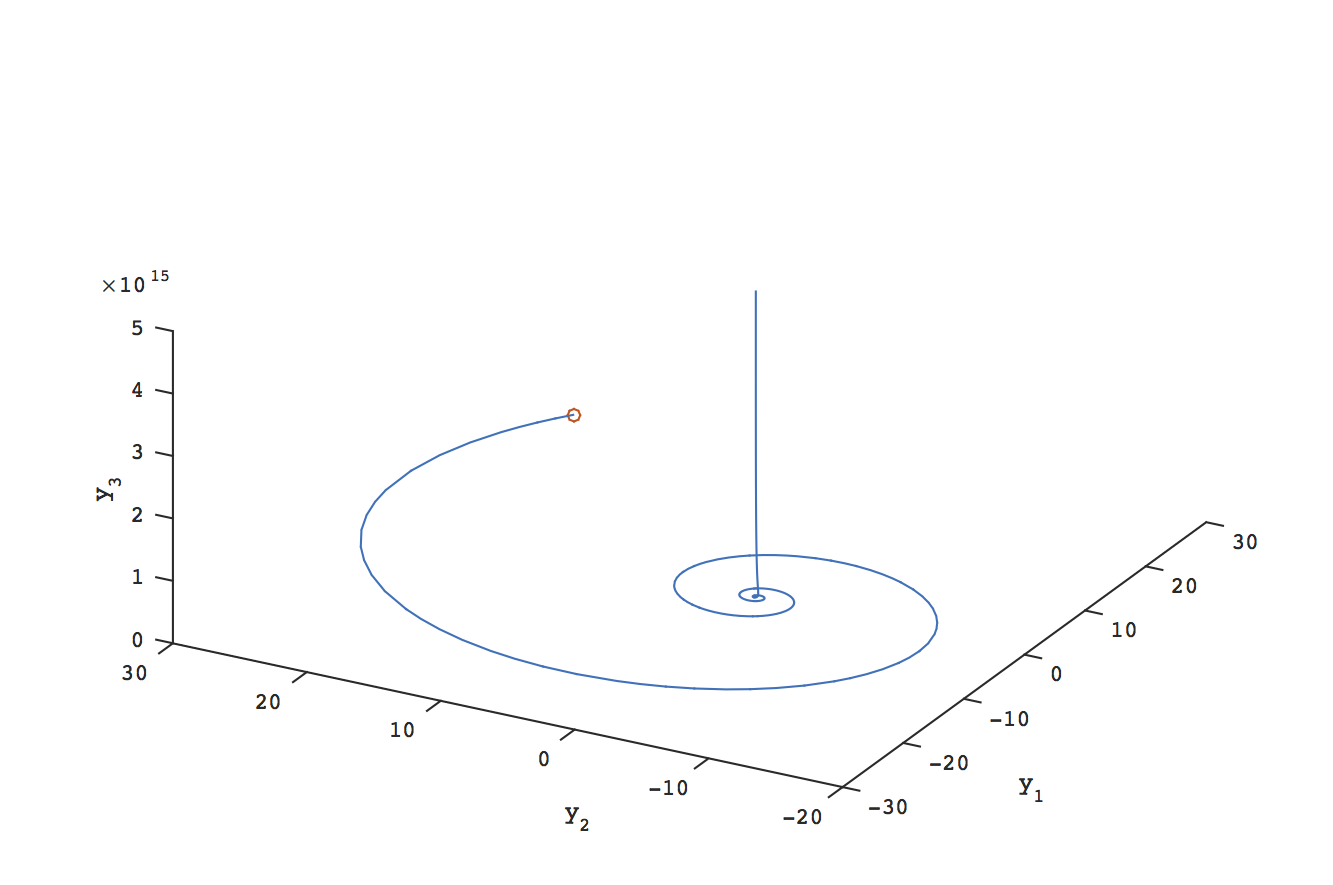}
\caption{Trajectory of a solution $y(t)$ of \eqref{eqn:ex3} when $a = -1.5$, $b = 1.0$, and $c = -1.25$.
}
\label{fig:behaviorex3}}
The initial point $x(0)=(25,25,25)$ is plotted by a circle marker. The trajectory first spirals up around the line $(0,0,y_3)$.
Then, it blows up when $y_1$ and $y_2$ become small.
\end{figure}

The third example is
\begin{equation}\label{eqn:ex3}
\left\{
\begin{array}{l}
	\displaystyle\frac{dy_1}{dt}=  (a-c)y_1y_3 - by_2y_3,\\[4mm]
	\displaystyle\frac{dy_2}{dt}= by_1y_3 + (a-c)y_2y_3,\\[4mm]
	\displaystyle\frac{dy_3}{dt}= -cy_3^2,
\end{array}
\right.
\end{equation}
where $a,b,c\not = 0$.
In this example, we also chose the Poincar\'e compactification in order to write the normalized vector field corresponding to \eqref{normalized-poincare}, which is given by
\begin{equation}\label{eqn:ex3_normalized}
\left\{
\begin{array}{l}
	\displaystyle\frac{dx_1}{d\tau}= (a-c)x_1x_3 - bx_2x_3 - G(x)x_1 ,\\[4mm]
	\displaystyle\frac{dx_2}{d\tau}= bx_1x_3 + (a-c)x_2x_3 - G(x)x_2,\\[4mm]
	\displaystyle\frac{dx_3}{d\tau}= -cx_3^2 - G(x)x_3,
\end{array}
\right.
\end{equation}
where $G(x) = ax_3(x_1^2+x_2^2) - cx_3(x_1^2+x_2^2+x_3^2)$.
Assume that $a, c < 0$, and $b\not = 0$. Then, the critical point at infinity, $(x_1,x_2,x_3) = (0,0,1)$, is a stable equilibrium.
One can see that the Jacobian matrix of the right-hand side of (\ref{eqn:ex3_normalized}) at $(0,0,1)$ is
\begin{equation*}
\begin{pmatrix}
a & -b & 0\\
b & a & 0\\
0 & 0 & 2c
\end{pmatrix}.
\end{equation*}
That is, the Jacobian matrix has one real negative eigenvalue and one complex conjugate pair of eigenvalues, and the critical point at infinity $(0,0,1)$ is a sink with a spiral.
Let $a = -1.5$, $b = 1.0$, and $c = -1.25$.
Figure \ref{fig:behaviorex3} illustrates the trajectory of \eqref{eqn:ex3}.

For $\varepsilon = 1.0\times 10^{-10}$, we set $N$ in \eqref{eqn:domainofLyapunov} using the critical point at infinity $x_\ast=(0,0,1)$.
We numerically verified that the Lyapunov function can be constructed on the set $\tilde{N}$ containing $N$.
Furthermore, we rigorously computed the solution of \eqref{eqn:ex3_normalized}, with the initial value $x(0)=(0.1,\,0.1,\,0.1)$ satisfying $x\to x_\ast$.
When $\tau_N=17$, we obtained
\[
{\small
x(\tau_N)\in
\left(
\begin{array}{c}
\left[6.5706769223025935\times 10^{-8}, 6.57067692230514\times 10^{-8}\right]\\
\left[-5.6195790918731807\times 10^{-8}, -5.619579091869991\times 10^{-8}\right]\\
\left[0.99999999995255339,0.99999999995256051\right]
\end{array}
\right)\subset N}
\]
and $t_N\in [7.9999713368258049,8.0000130782777568]$.
The Lyapunov function indicates that the solution $x$ with the initial value $x(0)$ asymptotically goes to the critical point at infinity $x_\ast$ as $\tau\to\infty$.
From \eqref{eqn:Blowuptime}, we have that $t_{\max}-t_N\le 5.5444899303784339\times 10^{-4}$.
Consequently, from Theorem \ref{th:main_poincare}, we have proved that the solution of \eqref{eqn:ex3} blows up, and its blow-up time is enclosed in the interval
\[
	t_{\max}\in[7.9999713368258049,8.0005675272707962].
\]
The execution time was about 0.778 seconds.
\subsection{Discretized nonlinear heat equation of order $3$}
\label{section-heat3}
The fourth example is the discretized heat equation with cubic nonlinearity:
\begin{equation}\label{eqn:ex4}
\left\{
\begin{array}{l}
	\displaystyle\frac{dy_1}{dt}= n^2(-2y_1 + y_2) + y_1^3,\\
	\qquad\vdots\\
	\displaystyle\frac{dy_k}{dt}= n^2(y_{k-1} - 2y_k + y_{k+1}) + y_{k}^3~~~(k=2,3,...,n-2),\\[2mm]
	\qquad\vdots\\
	\displaystyle\frac{dy_{n-1}}{dt}= n^2(y_{n-2} - 2y_{n-1}) + y_{n-1}^3,
\end{array}
\right.
\end{equation}
where $n$ is a positive integer.
This system, with the implicit boundary conditions $y_0=y_n = 0$, is considered as the finite difference semi-discretization of the nonlinear heat equation on a unit interval with the homogeneous Dirichlet boundary condition:
\[
\left\{
\begin{array}{l}
	u_t= u_{xx} + u^3,\quad t>0,\ x\in (0,1), \\
	u = 0  \text{ at }x=0,1.
\end{array}
\right.
\]
The grid points of the finite difference scheme are defined by $x_k=k/n$ $(k=0,1,...,n)$.
Then, it is expected that $y_k(t)$ approximates $u(t,x_k)$.

The normalized vector field associated with (\ref{eqn:ex4}) via the Poincar\'{e} compactification is
\begin{equation}\label{eqn:ex4_normalized}
\left\{
\begin{array}{l}
	\displaystyle\frac{dx_1}{d\tau}= \frac{n^2}{\kappa^{2}}(-2x_1 + x_2) + x_1^3 - G_3(x)x_1,\\
	\qquad\vdots\\
	\displaystyle\frac{dx_k}{d\tau}= \frac{n^2}{\kappa^{2}}(x_{k-1}-2x_k + x_{k+1}) + x_k^3 - G_3(x)x_k~~~(k=2,3,...,n-2),\\
	\qquad\vdots\\
	\displaystyle\frac{dx_{n-1}}{d\tau}= \frac{n^2}{\kappa^{2}}(x_{n-2}-2x_{n-1}) + x_{n-1}^3 - G_3(x)x_{n-1},
\end{array}
\right.
\end{equation}
where 
\begin{align*}
G_3(x) &= x_1\left(\frac{n^2}{\kappa^{2}}(-2x_1+x_2)+x_1^3\right) \\
	&\quad+ \sum_{k=2}^{n-2}x_k\left(\frac{n^2}{\kappa^{2}}(x_{k-1}-2x_k+x_{k+1})+x_k^3\right) + x_{n-1}\left(\frac{n^2}{\kappa^{2}}(x_{n-2}-2x_{n-1})+x_{n-1}^3\right),\\
	\kappa &= \kappa(x) = (1-\|x\|^2)^{-1/2}. 
\end{align*}

The initial value was set as $y_k(0)=10$ ($k=1,2,...,n-1$).
We rigorously computed the solution of \eqref{eqn:ex4_normalized} up to a fixed $\tau_N$.
The critical point at infinity $x_\ast$ was enclosed in an interval vector by using the Krawczyk method.
We also set $N$ with $\varepsilon$ such that $L\left(x(\tau_N)\right)<\varepsilon^2$.
Consequently, if the Lyapunov function is constructed on a set $\tilde N$ containing $N$, then it has been proved that the solution $x$ asymptotically goes to $x_\ast$ as $\tau\to\infty$.
On the basis of Theorem \ref{th:main_poincare}, we have obtained rigorous enclosures of the blow-up time, as given in Table \ref{tab:blowuptime_u3}.
\begin{table}[h]
\caption{Results of numerical validation for the blow-up time of \eqref{eqn:ex4}}
The blow-up time $t_{\max}$ is rigorously enclosed in the intervals given below.
Subscript and superscript numbers in the table denote lower and upper bound of the interval, respectively.
\begin{center}
\begin{tabular}{cccc}
\hline
$n$ & $\tau_N$ &inclusion of $t_{\max}$ & execution time\\
\hline\\[-2mm]
4 & 30 & $0.0050340400_{162383761}^{784869202}$ & 44.417s\\[2mm]
6 & 30 & $0.00500977_{0457049421}^{25547564119}$ & 3m15.793s\\[2mm]
8 & 30 & $0.005003_{7433760869625}^{9439361921953}$ & 7m11.506s\\[2mm]
10 & 35 & $0.005001_{7211768978893}^{9593060980734}$ & 20m56.604s\\[2mm]
12 & 40 & $0.00500_{08814990989457}^{12269722264098}$ & 72m10.282s\\[2mm]
14 & 50 & $0.005000_{4463436608779}^{6415371753669}$ & 202m0.581s\\[2mm]
16 & 50 & $0.00500_{0091354528135}^{12389304747536}$ & 387m26.199s\\[2mm]
\hline
\end{tabular}
\end{center}
\label{tab:blowuptime_u3}
\end{table}%
\subsection{The Riccati equation}\label{sec:riccati}
The fifth example again considers a simple ODE, the Riccati equation
\begin{equation}\label{eqn:ex5}
\frac{dy}{dt} = y^2 + t,~y(0)=\frac{1}{2},
\end{equation}
equivalently,
\[
\left\{
\begin{array}{l}
	\displaystyle\frac{dy_1}{dt} = y_1^2 + y_2, \\[4mm]
	\displaystyle\frac{dy_2}{dt} = 1
\end{array}
\right.,
\]
with the initial value $\left(y_1(0),\,y_2(0)\right)=(1/2,0)$.
When the Poincar\'{e} compactification is chosen, the normalized vector field of the form \eqref{normalized-poincare} associated with (\ref{eqn:ex5}) is given by
\[
\left\{
\begin{array}{l}
\displaystyle\frac{dx_1}{d\tau} = x_1^2 + (1-x_1^2-x_2^2)^{1/2}x_2 - G(x)x_1,\\[4mm]
\displaystyle\frac{dx_2}{d\tau} = (1-x_1^2-x_2^2) - G(x)x_2,
\end{array}
\right.
\]
where $G(x)=x_1^3 + (1-x_1^2-x_2^2)^{1/2} x_1x_2 + (1-x_1^2-x_2^2) x_2$.
Unlike in the previous examples, our validation method does not work in this case.
Indeed, the Jacobian matrix at $x = (x_1, x_2)$ is 
\begin{equation*}
\begin{pmatrix}
2x_1 - \frac{x_1x_2}{(1-x_1^2-x_2^2)^{1/2}} -G_{x_1}x_1-G(x) &  \frac{1-x_1^2 - 2x_2^2}{(1-x_1^2-x_2^2)^{1/2}} -G_{x_2}x_1 \\ -2x_1-G_{x_1}x_2 & -2x_2 - G_{x_2}x_2 - G(x)
\end{pmatrix},
\end{equation*}
where $G_{x_1}$ and $G_{x_2}$ denote the partial derivatives of $G(x)$ with respect to $x_1$ and $x_2$, respectively.
Because $\|x\| \to 1$ as $(x_1,x_2)\to \partial \mathcal{D}$, it turns out that the first row of the Jacobian matrix becomes infinite at $x_\ast$.
Thus, we cannot construct a Lyapunov function around the critical point at infinity.
\begin{rem}
\label{rem-failure}
The above failure follows mainly because the form of the nonlinearity $f(y)$ is 
\begin{equation*}
f(y) = p_0(y) + p_1(y) + \cdots + p_d(y),~p_d\not \equiv 0~\text{and}~p_{d-1}\not \equiv 0.
\end{equation*}
The term $p_{d-1}$ means that the coefficient $\kappa^{-1}$ in \eqref{normalized} and its differential may contain terms of order $O\left((1-\|x\|^2)^{-1/2}\right)$.
This causes a discontinuity of the Jacobian matrix to occur at $x_\ast$.
\end{rem}
To avoid singularities such as that above, we can apply the {\em parabolic} compactification.
See Section \ref{section-parabolic} for the details regarding parabolic compactifications.
In addition, see Section \ref{section-blowup-parabolic} for the validation of blow-up times.
The normalized vector field of the form \eqref{normalized-parabolic} associated with (\ref{eqn:ex5}) via the parabolic compactification is given by
\begin{equation}\label{eqn:ex5_normed}
\left\{
\begin{array}{l}
	\displaystyle\frac{dx_1}{d\tau} = (1+R^2)\{x_1^2+(1-R^2)x_2\} - 2\left \{x_1^4 + (1-R^2)x_1^2x_2 + (1-R^2)^2 x_1x_2\right \}, \\[4mm]
	\displaystyle\frac{dx_2}{d\tau} = (1+R^2)(1-R^2)^2 - 2\left \{ x_1^3x_2 + (1-R^2)x_1x_2^2 + (1-R^2)^2 x_2^2 \right \},
\end{array}
\right.
\end{equation}
where $R^2 = x_1^2+x_2^2$.
In this case, there is no failure in computing the Jacobian matrix at $x_\ast$.
The critical point at infinity is $x_\ast=(1,0)$.

For $\varepsilon = 1.0\times 10^{-10}$, we set $N$ as in \eqref{eqn:domainofLyapunov}.
We numerically verified that the Lyapunov function can be constructed on the set $\tilde{N}$ containing $N$.
When $\tau_N=7$, we obtained
\[
{\small x(\tau_N)\in
\left(
\begin{array}{c}
\left[0.99999653020284495,0.99999653020284563\right]\\
\left[9.9674655326613292\times 10^{-6},9.9674655326622847\times 10^{-6}\right]
\end{array}
\right)\subset N}
\]
and $t_N\in [1.4363412444327372,1.4363412444327977]$.
The Lyapunov function indicates that the solution $x$ of \eqref{eqn:ex5_normed} asymptotically goes to the critical point at infinity $x_\ast$ as $\tau\to\infty$.
From \eqref{eqn:Blowuptime_para}, we have that $t_{\max}-t_N\le 2.110827910648215\times 10^{-5}$.
Consequently, from Theorem \ref{th:main_para}, it has been proved that the solution of \eqref{eqn:ex5} blows up, and its blow-up time is enclosed in the interval
\[
	t_{\max}\in[1.4363412444327372,1.4363623527119043].
\]
The execution time was about 0.399 seconds.
\subsection{Discretized nonlinear heat equation of order $2$}
The final example is the discretized heat equation with second-order nonlinearity:
\begin{equation}\label{eqn:ex6}
\left\{
\begin{array}{l}
	\displaystyle\frac{dy_1}{dt}= n^2(-2y_1 + y_2) + y_1^2,\\
	\qquad\vdots\\
	\displaystyle\frac{dy_k}{dt}= n^2(y_{k-1} - 2y_k + y_{k+1}) + y_{k}^2~~~(k=2,3,...,n-2),\\[2mm]
	\qquad\vdots\\
	\displaystyle\frac{dy_{n-1}}{dt}= n^2(y_{n-2} - 2y_{n-1}) + y_{n-1}^2,
\end{array}
\right.,
\end{equation}
where $n$ is a positive integer.
This system, with the implicit boundary conditions $y_0=y_n = 0$, is considered as 
the spatially finite difference discretization of the nonlinear heat equation on the unit interval with the homogeneous Dirichlet boundary condition
\begin{equation*}
\left\{
\begin{array}{l}
	u_t= u_{xx} + u^2,\quad t>0,\ x\in (0,1), \\
	u = 0 \text{ at }x=0,1.
\end{array}
\right.
\end{equation*}
The grid points of the finite difference scheme are defined by $x_k=k/n$ $(k=0,1,...,n)$.
Then, it is expected that $y_k(t)$ in \eqref{eqn:ex6} approximates $u(t,x_k)$.

Because this system falls under the case described in Remark \ref{rem-failure} with $d=2$, the Poincar\'{e} compactification is not available for validating blow-up solutions.
Thus, we apply the parabolic compactification.
The normalized vector field associated with (\ref{eqn:ex6}) via the parabolic compactification is given by
\begin{equation}\label{eqn:ex6_normalized}
\left\{
\begin{array}{l}
	\displaystyle\frac{dx_1}{d\tau}= \left(1+R^2\right)\left\{n^2\left(1-R^2\right)\left(-2x_1 + x_2\right) + x_1^2\right\} - 2G_2(x)x_1,\\
	\qquad\vdots\\
	\displaystyle\frac{dx_k}{d\tau}= \left(1+R^2\right)\left\{n^2\left(1-R^2\right)\left(x_{k-1}-2x_k + x_{k+1}\right) + x_k^2\right\} - 2G_2(x)x_k~~~(k=2,3,...,n-2),\\
	\qquad\vdots\\
	\displaystyle\frac{dx_{n-1}}{d\tau}= \left(1+R^2\right)\left(n^2\left(1-R^2\right)\left(x_{n-2}-2x_{n-1}\right) + x_{n-1}^2\right) - 2G_2(x)x_{n-1},
\end{array}
\right.
\end{equation}
where 
\begin{align*}
G_2(x) &= x_1\left\{n^2\left(1-R^2\right)\left(-2x_1+x_2\right)+x_1^2\right\}
	+ \sum_{k=2}^{n-2}x_k\left\{n^2\left(1-R^2\right)\left(x_{k-1}-2x_k+x_{k+1}\right)+x_k^2\right\}\\
	&\quad+ x_{n-1}\left\{n^2\left(1-R^2\right)\left(x_{n-2}-2x_{n-1}\right)+x_{n-1}^2\right\},\\
	R^2 &= \sum_{j=1}^{n-1} x_j^2.
\end{align*}

The initial value was set as $y_k(0)=10$ ($k=1,2,...,n-1$).
By using validated computations, we rigorously enclosed the solution of \eqref{eqn:ex6_normalized} up to a fixed $\tau_N$.
The critical point at infinity, $x_\ast$, was enclosed in an interval vector using the Krawczyk method.
We also set $N$ with $\varepsilon$ such that $L\left(x(\tau_N)\right)<\varepsilon^2$.
Consequently, if the Lyapunov function is constructed on a set $\tilde N$ containing $N$, then it has been proved that the solution $x$ asymptotically goes to $x_\ast$ as $\tau\to\infty$.
From Theorem \ref{th:main_para}, we have obtained a rigorous enclosure of the blow-up time, which is presented in Table \ref{tab:blowuptime_u2}.
{
We note that, if we consider the blow-up time of \eqref{eqn:ex6} (and \eqref{eqn:ex4}) for large $n$, it is difficult to find $N$ that becomes the Lyapunov domain defined in Proposition \ref{prop:Lyapunov}.
The Lyapunov domain is validated if the negative definiteness of the matrix $A(x)$ for all $x\in N$ defined in \eqref{eqn:ax} is ensured.
If $N$ is large, i.e., $\varepsilon$ is large, then the validation of the negative definiteness is expected to fail.
This is the reason of failure in our validating method when $n=16$ with $\tau_N=20$.
}

\begin{table}[h]
\caption{Results of numerical validation for the blow-up time of \eqref{eqn:ex6}}
The blow-up time $t_{\max}$ is rigorously enclosed in the intervals given below.
Subscript and superscript numbers in the table denote the lower and upper bounds of the interval, respectively.
When $n=16$ with $\tau_N=20$, we could not validate any Lyapunov domain (defined in Proposition \ref{prop:Lyapunov}) around $x_\ast$.
Therefore, we concluded that our validating method failed to enclose the blow-up time of \eqref{eqn:ex6}.
\begin{center}
\begin{tabular}{cccc}
\hline
$n$ & $\tau_N$ & inclusion of $t_{\max}$ & execution time\\
\hline\\[-2mm]
4 & 15 & $0.242_{86876161046069}^{90697501550363}$ & 32.947s\\[2mm]
6 & 20 & $0.2462_{3855107071979}^{4064886491729}$ & 5m3.502s\\[2mm]
8 & 20 & $0.24608_{006592024286}^{664310433196}$ & 15m2.741s\\[2mm]
10 & 20 & $0.245_{78076022169118}^{80239357994319}$ & 53m35.880s\\[2mm]
12 & 20 & $0.245_{55282235874756}^{6169375007353}$ & 126m22.870s\\[2mm]
14 & 20 & $0.245_{39149499338358}^{5655918000558}$ & 273m48.783s\\[2mm]
16 & 20 & Failed & 492m54.807s\\[2mm]
\hline
\end{tabular}
\end{center}
\label{tab:blowuptime_u2}
\end{table}%

{
\begin{rem}
When all componets of the critical point at infinity $x_\ast$ are nonzero, we can say that all the components of the solution $y$ blow up at the same time.
This is because the scaling is chosen to be uniform in the sense that all the components grow in the same rate, e.g., as seen in the definition of the Poincar\'e compactification, say
\begin{equation*}
y_i = \kappa\left(T^{-1}(x)\right)x_i=\frac{x_i}{(1-\|x\|^2)^{1/2}}\quad \text{ for all }i.
\end{equation*}
Furthermore, the definition of critical points at infinity implies that 
\begin{equation*}
|y_i| \xrightarrow[\tau \to\infty]{} \infty,\quad \frac{y_i}{\|y\|}\xrightarrow[\tau \to \infty]{} (x_\ast)_i \not = 0.
\end{equation*}
Since $\tau \to \infty$ corresponds to $t \to t_{\max}$ if our validation succeeds, the growth rate of all $y_i$ is the same as $\|y\|$ in the $\tau$-time scale.
On the other hand, there is a different situation if (at least) one of the components of $x_\ast$ is zero, say $(x_\ast)_i = 0$.
In this case, the growth rate of $y_i$ is smaller than that of $\|y\|$.
We can thus say that $y_i$ grows up in the smaller rate than $\|y\|$ (in the $\tau$-time scale), which is finite or may go to infinity, as indicated in Section \ref{sec:riccati}.
\end{rem}
}

\section*{Conclusion}
In this paper, we have discussed a method for validating blow-up solutions of ordinary differential equations.
Our method consists of a geometric transformation of phase spaces, called a {\em compactification}; the use of {\em Lyapunov functions}; and the standard integration of differential equations.
These concepts explicitly represent points and dynamics at ``infinity'' as computable ones, which enable us to validate
\begin{itemize}
\item blow-up solutions, and
\item lower and upper bound on their blow-up times,
\end{itemize}
with computer assistance. 
Unlike in previous numerical studies regarding blow-up solutions, our method {\em rigorously} validates blow-up solutions, and encloses blow-up times.
We have also provided several numerical validation results, to illustrate the applicability of our method, involving not only autonomous but also non-autonomous systems.
We believe that our proposed approach provides new insight into the analysis and the computation of blow-up solutions of differential equations.

\bigskip
We conclude this paper by commenting on the possible extension of our method in a particular direction.
A natural question is whether our method can be extended to blow-up solutions of {\em partial differential equations}.
Preceding works indicate that compactification itself makes sense for Hilbert spaces (cf. \cite{bib:H}), and that the Lyapunov functions around equilibria can be validated for parabolic PDEs (cf. \cite{bib:Mat}).
Therefore, it is natural to expect that our procedure can be extended to PDEs, at least to parabolic PDEs.
However, several difficulties arise regarding such an extension.

The first problem is the detection of isolated critical points at infinity, which is essential for constructing the Lyapunov functions.
Even for polynomial PDEs, such as the nonlinear heat equation $u_t = u_{xx} + u^d$ in the unit interval $[0,1]\subset \mathbb{R}$ with the homogeneous Dirichlet boundary condition, direct calculations of the normalized vector field via the Fourier expansion $u(t,x) = \sum_{k=1}^\infty a_k(t)\sin(k\pi x)$, for example in \cite{bib:H}, imply that there are infinitely many critical points at infinity.
This situation implies that critical points at infinity may accumulate.

The second problem is the direct calculation of trajectories that blow up using validated computations.
Several methodologies exist for computing rigorous trajectories of PDEs, such as in \cite{bib:AK2010, bib:CGL2015, bib:CZ2015, bib:MTKO2015, bib:Zgl2010},
all of which focus on validations of {\em bounded} trajectories.
All of these methodologies employ either power estimates of eigenvalues or the properties of analytic semigroups. 
In both cases, the ``decay'' property stemming from the leading linear operators is essential for the validations.
On the other hand, nonlinear growth becomes dominant for the time evolution of blow-up solutions (see, e.g., \cite{bib:FM2002, bib:Meier1990}).
This fact presents the non-trivial question of whether the preceding validation techniques for PDEs can be really applied to blow-up solutions of PDEs.
\section*{Acknowledgments.}
AT was partially supported by JSPS Grant-in-Aid for Young Scientists (B), No.\,15K17596.
KM was partially supported by The Coop with Math Program, a commissioned project by Ministry of Education, Culture, Sports, Science and Technology (MEXT).
SO was partially supported by CREST, JST.
The authors also wish to thank to Prof. Nobito Yamamoto, who provided us with the essential suggestion of re-parameterization in the integrals of ODEs via Lyapunov functions.

\bibliographystyle{plain}
\bibliography{blow_up}


\end{document}